\documentclass[final]{siamltex}

\usepackage{amsmath} 
\usepackage{amsfonts}

\usepackage{graphicx}
\usepackage{subcaption}
\usepackage{url}
\usepackage{color}

\newtheorem{remark}{Remark}
\newtheorem{assumption}{Assumption}

\title{Fast Iterative Solver for the Optimal Control of Time-Dependent PDEs with Crank--Nicolson Discretization in Time}

\author{Santolo Leveque\thanks{School of Mathematics, The University of Edinburgh, James Clerk Maxwell Building, The King's Buildings, Peter Guthrie Tait Road, Edinburgh, EH9 3FD, United Kingdom ({\tt S.Leveque@sms.ed.ac.uk})} \and John W. Pearson\thanks{School of Mathematics, The University of Edinburgh, James Clerk Maxwell Building, The King's Buildings, Peter Guthrie Tait Road, Edinburgh, EH9 3FD, United Kingdom ({\tt j.pearson@ed.ac.uk})}}

\begin{document}
\maketitle

\begin{abstract}
In this article, we derive a new, fast, and robust preconditioned iterative solution strategy for the all-at-once solution of optimal control problems with time-dependent PDEs as constraints, including the heat equation and the non-steady convection--diffusion equation. After applying an optimize-then-discretize approach, one is faced with continuous first-order optimality conditions consisting of a coupled system of PDEs. As opposed to most work in preconditioning the resulting discretized systems, where a (first-order accurate) backward Euler method is used for the discretization of the time derivative, we employ a (second-order accurate) Crank--Nicolson method in time. We apply a carefully tailored invertible transformation for symmetrizing the matrix, and then derive an optimal preconditioner for the saddle-point system obtained. The key components of this preconditioner are an accurate mass matrix approximation, a good approximation of the Schur complement, and an appropriate multigrid process to apply this latter approximation---these are constructed using our work in transforming the matrix system. We prove the optimality of the approximation of the Schur complement through bounds on the eigenvalues, and test our solver against a widely-used preconditioner for the linear system arising from a backward Euler discretization. These demonstrate the effectiveness and robustness of our solver with respect to mesh-sizes, regularization parameter, and diffusion coefficient.
\end{abstract}



\pagestyle{myheadings}
\thispagestyle{plain}
\markboth{S. LEVEQUE AND J. W. PEARSON}{ITERATIVE SOLVER FOR OPTIMAL CONTROL OF TIME-DEPENDENT PDES}

\section{Introduction}
\label{sec1}

PDE-constrained optimization is a subject area which has attracted significant interest in recent years (see \cite{Hinze,Troltzsch} for a comprehensive overview of the field). Due the complex structure and high dimensionality of these problems when an accurate discrete solution is sought, preconditioned iterative solvers have been employed for the `all-at-once' resolution of such formulations, see for example \cite{Rees_Dollar_Wathen, Schoberl_Zulehner} for early work in this direction. In this paper, we consider fast and robust preconditioned solvers for matrix systems arising from the \emph{optimize-then-discretize} approach applied to time-dependent PDE-constrained optimization problems. We examine the optimal control of the heat equation with Dirichlet boundary conditions, to allow benchmarking of our method against a widely-applied preconditioned iterative method for its solution, and time-dependent convection--diffusion control problems.

When preconditioners are sought for certain time-dependent problems, it is typical to apply a (first-order accurate) backward Euler method in time, as this leads to particularly convenient structures within the matrix and facilitates effective preconditioning; see \cite{Pearson_Stoll_Wathen} for a mesh- and $\beta$-robust preconditioner for the heat control problem, and \cite{Dolgov_Stoll, Pearson_Stoll, Stoll_Benner_Onwunta_Dolgov, Stoll_Pearson_Maini, Yucel_Stoll_Benner} for applications to different problems. However, the required discretization strategy results in slower convergence in time than space: if a method is second-order accurate in space, it is reasonable to choose $\tau = \mathcal{O}(h^2)$, where $\tau$ is the time step and $h$ is the mesh-size in space. This results in a matrix system of huge dimension, and hence a high CPU time is required for its solution. The question we wish to investigate here is whether applying a higher-order Crank--Nicolson method in time is beneficial, due to the reduced dimension of the matrix system required to obtain a similar discretization error. The challenge here is that the much more complex structure of the resulting matrix system makes preconditioning a highly non-trivial task, so a more sophisticated numerical method and preconditioning strategy need to be devised to achieve fast and robust convergence of the iterative solver.

The remainder of this article is structured as follows. In Section \ref{2} we introduce the problems we consider, that is time-dependent convection--diffusion control and heat control, and outline the matrices arising upon discretization of such problems as well as a stabilization technique used for convection--diffusion control problems. In Section \ref{3} we describe the matrix systems obtained upon discretization of the first-order optimality conditions, and in particular demonstrate a suitable transformation which allows symmetrization of the matrix system obtained from the Crank--Nicolson method. This enables us to apply a symmetric iterative solver such as MINRES \cite{Paige_Saunders}, which is highly desirable from the perspective of proving convergence of the iterative method. In Section \ref{4} we derive our proposed new preconditioner, using saddle-point theory along with suitable approximations of the $(1,1)$-block and Schur complement, and provide eigenvalue results for the preconditioned matrix system. In Section \ref{5} we benchmark our method against the widely-used backward Euler method, coupled with the preconditioner derived in \cite{Pearson_Stoll_Wathen}, and demonstrate our new preconditioner's efficiency and robustness with respect to all parameters involved in the convection--diffusion control problem.

\section{Problem Formulation}
\label{2}
In this article we consider the fast and robust numerical solution of time-dependent PDE-constrained optimization problems. In particular, we examine distributed convection--diffusion control problems of the form:
\begin{equation}\label{convection_diffusion_control_functional}
\min_{y,u}~~J(y,u) = \dfrac{1}{2} \int_{0}^{t_{f}} \int_{\Omega} |y(\mathbf{x},t)- \widehat{y}(\mathbf{x},t)|^2 \: {\rm d}\Omega {\rm d}t + \dfrac{\beta}{2} \int_{0}^{t_{f}} \int_{\Omega} |u(\mathbf{x},t)|^2 \: {\rm d}\Omega {\rm d}t
\end{equation}
subject to
\begin{equation}\label{convection_diffusion_control_constraints}
\left\{
\begin{array}{rl}
\vspace{1ex}
\dfrac{\partial y}{\partial t} - \epsilon \nabla^2 y + \mathbf{w} \cdot \nabla y = u + f(\mathbf{x},t) & \quad \mathrm{in} \; \Omega \times (0,t_{f}), \\
\vspace{1ex}
y(\mathbf{x},t) = g_D(\mathbf{x},t) & \quad \mathrm{on} \; \partial \Omega_D \times (0,t_{f}),\\
\vspace{1ex}
\dfrac{\partial y}{\partial n}(\mathbf{x},t) = g_N(\mathbf{x},t) & \quad \mathrm{on} \; \partial \Omega_N \times (0,t_{f}),\\

y(\mathbf{x},0) = y_0(\mathbf{x}) & \quad \mathrm{in} \; \Omega,

\end{array}
\right.
\end{equation}
where the variables $y$, $\widehat{y}$, and $u$ are the \emph{state}, \emph{desired state}, and \emph{control} variables, respectively, $\beta>0$ is a \emph{regularization parameter}, $\epsilon>0$ is the diffusion coefficient, and $\mathbf{w}$ is a divergence-free wind vector (i.e., $\nabla \cdot \mathbf{w} = 0$). The problem is solved on a spatial domain $\Omega \subset \mathbb{R}^{d}$, $d\in\{1,2,3\}$, with boundary $\partial \Omega = \partial \Omega_D \cup \partial \Omega_N$, $\partial \Omega_D \cap \partial \Omega_N = \emptyset$, up to a final time $t_{f}>0$, that is $(\mathbf{x},t)\in\Omega \times (0,t_{f})$. Here $\dfrac{\partial y}{\partial n}(\mathbf{x},t)$ represents the normal derivative of $y$ on $\partial \Omega_N$. The functions $f, \: g_D, \: g_N$, and $y_0$ are known.

In \eqref{convection_diffusion_control_constraints}, the term $- \epsilon \nabla^2 y$ denotes the diffusive element, and the term $\mathbf{w} \cdot \nabla y$ represents convection. In physical (real-world) problems, as pointed out for example in \cite{Elman_Silvester_Wathen}, convection typically plays a more significant physical role than diffusion, so $\epsilon \ll \|\mathbf{w} \|$ for many practical problems. However this in turn makes the problem more difficult to solve \cite{Elman_Silvester_Wathen, Ramage} as the solution procedure will need to be robust with respect to the direction of the wind $\mathbf{w}$ and any boundary or internal layers that form. The presence of boundary or internal layers is also an issue that we have to deal with when discretizing the convection--diffusion differential operator. Indeed, when solving a convection--diffusion problem, a stabilization procedure is often utilized in order to `capture' all the layers.

Though we will test our solver on the convection--diffusion control problem \eqref{convection_diffusion_control_functional}--\eqref{convection_diffusion_control_constraints}, a particular case of interest arises when $\epsilon=1$, $\mathbf{w}=\mathbf{0}$, and $\partial \Omega_D=\partial \Omega, \: \partial \Omega_N = \emptyset$, giving the following heat control problem with Dirichlet boundary conditions:
\begin{align}
\label{heat_control_functional} \min_{y,u}&\quad{}J(y,u) \\
\label{heat_control_constraints} \text{s.t.}&\quad\left\{
\begin{array}{rl}
\vspace{1ex}
\dfrac{\partial y}{\partial t} - \nabla^2 y  = u + f(\mathbf{x},t) & \quad \mathrm{in} \; \Omega \times (0,t_{f}), \\
\vspace{1ex}
y(\mathbf{x},t) = g(\mathbf{x},t) & \quad \mathrm{on} \; \partial \Omega \times (0,t_{f}),\\
y(\mathbf{x},0) = y_0(\mathbf{x}) & \quad \mathrm{in} \; \Omega.
\end{array}
\right.
\end{align}
This problem provides a valuable benchmark of our method against widely-used preconditioned iterative solution with a backward Euler method in time \cite{Pearson_Stoll_Wathen}. We note that our method can be readily generalized to heat control problems with Neumann and mixed boundary conditions.

\subsection{Discretization matrices and stabilization}
\label{2_1}

To illustrate the matrices involved in the finite element discretizations of optimal control problems, consider a standard Galerkin finite element discretization for the (steady-state) convection--diffusion problem:
\begin{equation}\label{forward_convection_diffusion_control_constraints}
- \epsilon \nabla^2 y + \mathbf{w} \cdot \nabla y = u + f(\mathbf{x}) \qquad \mathrm{in} \; \Omega \:.
\end{equation}
Letting $\{\phi_i\}_{i=1}^{n_x}$ be the same finite element basis functions for $y$ and $u$, then we would like to find approximations $y(\mathbf{x}) \approx \sum_{i=1}^{n_x} y_i \phi_i$, $u(\mathbf{x}) \approx \sum_{i=1}^{n_x} u_i \phi_i$. Letting the vectors $\mathbf{y} = \{y_i\}_{i=1}^{n_x}$, $\mathbf{u} = \{u_i\}_{i=1}^{n_x}$, a discretized version of \eqref{forward_convection_diffusion_control_constraints} is
\begin{equation*}
L \mathbf{y} := (\epsilon K + N + W) \mathbf{y} = M \mathbf{u} + \mathbf{f},
\end{equation*}
where (excluding the effect of boundary conditions):
\begin{align*}
K=\{k_{ij}\}_{i,j=1}^{n_x}, \qquad & k_{ij} = \int_{\Omega} \nabla \phi_i \cdot \nabla \phi_j \: \mathrm{d} \Omega \:, \\
N=\{n_{ij}\}_{i,j=1}^{n_x}, \qquad & n_{ij} = \int_{\Omega} ( \mathbf{w} \cdot \nabla \phi_j) \phi_i \: \mathrm{d} \Omega \:, \\
M=\{m_{ij}\}_{i,j=1}^{n_x}, \qquad & m_{ij} = \int_{\Omega} \phi_i \phi_j \: \mathrm{d} \Omega \:, \\
\mathbf{f}=\{f_{i}\}_{i=1}^{n_x}, \qquad & f_i = \int_{\Omega} f \phi_i \: \mathrm{d} \Omega \:,
\end{align*}
and the matrix $W$ denotes a possible stabilization matrix for the convection operator. For a heat control problem of the form \eqref{heat_control_constraints}, $\epsilon=1$, and $N=W=0$. We note that $K$ is generally referred to as a \emph{stiffness matrix}, and is symmetric positive semi-definite (positive definite unless $\partial\Omega_D = \emptyset$), and $M$ is referred to as a \emph{mass matrix}, which is symmetric positive definite. The matrix $N$ is skew-symmetric (meaning $N+N^\top=0$) in the case of Dirichlet problems; otherwise we obtain that
\begin{equation}\label{matrix_H}
n_{ij} + n_{ji} = \int_{\partial \Omega_N} \phi_i \phi_j \mathbf{w} \cdot \mathbf{n} \:\mathrm{d} s \:,
\end{equation}
using Green's theorem. In this latter case the spectral properties of the matrix $H:=N+N^\top$, which can be indefinite, will be useful for our subsequent analysis. We note that $\mathbf{w}\cdot \mathbf{n} = \left\| \mathbf{w} \right\| \cos{\theta}$, where $\theta$ denotes the angle between the wind $\mathbf{w}$ and the (outer) unit normal vector $\mathbf{n}$ on the boundary $\partial \Omega_N$, so setting $c=\max_{\mathbf{x}}\left\| \mathbf{w} \right\|$, we have that
\begin{equation}\label{wind_normal}
- c \leq \mathbf{w} \cdot \mathbf{n} \leq c.
\end{equation}
Defining the matrix
\begin{displaymath}
M_{\partial \Omega_N} = \{m_{ij}^{\partial \Omega_N}\}_{i,j=1}^{n_x}, \qquad m_{ij}^{\partial \Omega_N} = \int_{\partial \Omega_N} \phi_i \phi_j \:\mathrm{d} s \:,
\end{displaymath}
we may then write
\begin{displaymath}
\begin{array}{ll}
\vspace{1ex}
\mathbf{v}^\top (H + c M_{\partial \Omega_N}) \mathbf{v} & = \displaystyle \sum_{i=1}^{n_x} \sum_{j=1}^{n_x} v_i \left[ \int_{\partial \Omega_N} \left( \mathbf{w} \cdot \mathbf{n} + c \right) \phi_i \phi_j  \: \mathrm{d} s \right] v_j \\
& = \displaystyle \int_{\partial \Omega_N} \left( \mathbf{w} \cdot \mathbf{n} + c \right) v^2 \: \mathrm{d} s \geq 0
\end{array}
\end{displaymath}
for any $\mathbf{0} \neq \mathbf{v}\in \mathbb{R}^{n_x}$, where we set $v:= \sum_{i=1}^{n_x} v_i \phi_i$ and use \eqref{wind_normal}. Therefore,
\begin{equation}\label{matrix_H_mass_matrix}
H \succeq - c M_{\partial \Omega_N},
\end{equation}
where the notation $\Psi_1 \succeq \Psi_2$ means $\Psi_1 - \Psi_2$ is positive semi-definite. Hence, any negative eigenvalues of $H$ are dominated by a constant factor of those of the mass matrix related to the basis functions on $\partial \Omega_N$. This observation will be useful when discussing our approach for problems with Neumann or mixed boundary conditions.

Concerning the stabilization scheme used for solving the forward convection--diffusion problem, a popular stabilized finite element method is the Streamline Upwind Petrov--Galerkin (SUPG) method \cite{Hughes_Brooks}. However, as pointed out in \cite{Collis_Heinkenschloss, Heinkenschloss_Leykekhman}, applying this scheme to the (steady-state) control problem gives rise to extra difficulties. Specifically, applying the scheme to the control problem with the \textit{discretize-then-optimize} strategy leads to symmetric discrete equations in which the discrete adjoint problem is not a consistent discretization of the continuous adjoint problem, whereas the optimize-then-discretize approach gives rise to a different, non-symmetric discretized system which therefore does not possess the structure of a discrete optimization problem. In this work, we thus employ the adjoint-consistent Local Projection Stabilization (LPS) approach described in \cite{Becker_Braack, Becker_Vexler, Braack_Burman,Pearson_Wathen}, for which the discretization and optimization steps commute in the stationary case. Further, this stabilized finite elements leads to an order of convergence of $\mathcal{O}(h^{3/2})$ for the $L^2$-error \cite{Becker_Vexler}, which is optimal on general quasi-uniform meshes for the forward problem, see e.g., \cite{Zhou}. In the LPS formulation, the stabilization matrix $W$ is defined as
\begin{equation}\label{LPS_stab_matrix}
W=\{w_{ij}\}_{i,j=1}^{n_x}, \qquad w_{ij} = \delta \int_\Omega [\mathbf{w} \cdot \nabla \phi_i - \pi_h(\mathbf{w} \cdot \nabla \phi_i)] [\mathbf{w} \cdot \nabla \phi_j - \pi_h(\mathbf{w} \cdot \nabla \phi_j)] \: \mathrm{d} \Omega \:.
\end{equation}
Here, $\delta>0$ denotes a stabilization parameter, and $\pi_h$ is an orthogonal projection operator. From \eqref{LPS_stab_matrix}, the matrix $W$ can be viewed as a shifted discrete diffusion operator associated with the
streamline direction defined by $\mathbf{w}$. It is symmetric and positive semi-definite, as shown by following the working in \cite[p.\,17]{Elman_Silvester_Wathen}: letting $\mathbf{0}\neq\mathbf{v}\in \mathbb{R}^{n_x}$, and setting $\pi_h^i = \pi_h(\mathbf{w} \cdot \nabla \phi_i)$, $v := \sum_{i=1}^{n_x} v_i \phi_i$, $\widetilde{\pi}:= \sum_{i=1}^{n_x} v_i \pi_h^i$, we have
\begin{align*}
\mathbf{v}^\top W \mathbf{v} ={}& \delta \sum_{i=1}^{n_x} \sum_{j=1}^{n_x} v_i \left[ \int_\Omega [\mathbf{w} \cdot \nabla \phi_i - \pi_h^i] [\mathbf{w} \cdot \nabla \phi_j - \pi_h^{\:j}] \: \mathrm{d} \Omega \right] v_j\\
 ={}& \delta \int_{\Omega} \left(\mathbf{w} \cdot \nabla v - \widetilde{\pi}\right)\left(\mathbf{w} \cdot \nabla v - \widetilde{\pi}\right) \mathrm{d} \Omega = \delta \, \| \mathbf{w} \cdot \nabla v - \widetilde{\pi} \|_{L^2(\Omega)}^2 \geq 0,
\end{align*}
where we have used that $\sum_{i=1}^{n_x} v_i \nabla \phi_i = \nabla v$.

For the convergence of the method, we require $\pi_h$ to be an $L^2$-orthogonal (discontinuous) projection operator defined on patches of the domain $\Omega$ that satisfies the approximation and stability properties specified in \cite{Becker_Vexler}, where by a patch we mean the union of elements of our finite element discretization; the projection operator $\pi_h$ is left free to be discontinuous on the edges of the patches. In our implementation we will make use of $\mathbf{Q}_1$ elements, so the domain is divided into patches consisting of 2 elements in each dimension. To ensure the aforementioned properties are satisfied, as in \cite{Becker_Braack} we define $\pi_h$ as
\begin{displaymath}
\left.\pi_h(v)\right|_{\mathtt{P}} = \dfrac{1}{|\mathtt{P}|} \int_{\mathtt{P}} v \: \mathrm{d} \mathtt{P}, \quad \forall v \in L^2(\Omega),
\end{displaymath}
where $\left.\pi_h(v)\right|_{\mathtt{P}}$ is the restriction of $\pi_h(v)$ to the patch $\mathtt{P}$, and $|\mathtt{P}|$ is the (Lebesgue) measure of the patch. We refer again to \cite{Becker_Vexler} for the theoretical proof of the convergence of this method with this definition of the local projection operator. As in \cite[p.\,253]{Elman_Silvester_Wathen}, we choose $\delta$ locally on each patch $\mathtt{P}_k$ as $\delta_k$, with
\begin{displaymath}
\delta_k =
\left\{
\begin{array}{lcl}
\vspace{1ex}
\dfrac{h_k}{2 \| \mathbf{w}_k \| } \left(1 - \dfrac{1}{P_k}\right)  & \quad \mathrm{if} \: P_k > 1, \\
0 & \quad \mathrm{if} \: P_k \leq 1,
\end{array}
\right.
\end{displaymath}
where $\| \mathbf{w}_k \|$ is the $\ell^2$-norm of the wind at the patch centroid, $h_k$ is a measure of the patch length in the direction of the wind, and $P_k=\frac{\| \mathbf{w}_k \| h_k}{2 \epsilon}$ is the \textit{patch P\'{e}clet number}.

We observe that, with this (non-constant) choice of $\delta$, the matrix $W$ is still positive semi-definite. To prove this, it is sufficient to define $\pi_h^{i,k} = \pi_h(\mathbf{w} \cdot \nabla \phi_i)|_{{\mathtt{P}_k}}$ and $\widetilde{\pi}_k:= \sum_{i=1}^{n_x} v_i \pi_h^{i,k}$ locally on each patch and then proceed as above, obtaining
\begin{displaymath}
\displaystyle \mathbf{v}^\top W \mathbf{v} = \sum_k \delta_k  \| \mathbf{w} \cdot \nabla v - \widetilde{\pi}_k \|_{L^2(\mathtt{P}_k)}^2 \geq 0.
\end{displaymath}
The spectral properties of the matrices $K$, $M$, and $W$ (whether Dirichlet, Neumann or mixed boundary conditions are imposed) will be useful later, when discussing the optimality of our preconditioning approach.

Informed by the definitions of this section, we make the following assumption when carrying out the theoretical analysis in the remainder of this paper, and later discuss how our methodology could be applied if the assumption is relaxed:
\begin{assumption}\label{assumption_w}
When a convection-diffusion control problem is being solved we assume $\mathbf{w}$ is such that $\mathbf{w}\cdot\mathbf{n}=0$ on $\partial\Omega_N$. We allow heat control problems to satisfy any boundary conditions of the form \eqref{convection_diffusion_control_constraints}, and for convection--diffusion control we remove this restriction on the wind $\mathbf{w}$ for a pure Dirichlet problem (i.e., $\partial\Omega_N=\emptyset$).
\end{assumption}

When Assumption \ref{assumption_w} holds, the matrix $H=0$. There are a number of wind vectors $\mathbf{w}$ that satisfy the property above on the whole of $\partial\Omega$, see for example the well-known `recirculating wind' example described in \cite[p.\,240]{Elman_Silvester_Wathen}.

\section{First-Order Optimality Conditions and Discretizations in Time}
\label{3}

We now describe the strategy used for obtaining an approximate solution of \eqref{convection_diffusion_control_functional}--\eqref{convection_diffusion_control_constraints}. We apply an all-at-once approach coupled with the optimize-then-discretize scheme, in which the continuous Lagrangian is used to arrive at first-order optimality conditions, which are then discretized. For simplicity the working of this section considers Dirichlet boundary conditions, that is $\partial \Omega_D=\partial \Omega$, but it may be readily extended to problems where $\partial \Omega_D \subset \partial \Omega$. Introducing the \emph{adjoint} variable (or \emph{Lagrange multiplier}) $p$, we consider the Lagrangian associated to \eqref{convection_diffusion_control_functional}--\eqref{convection_diffusion_control_constraints} as in \cite{Becker_Vexler}. Then, by deriving the Karush--Kuhn--Tucker (KKT) conditions, the solution of \eqref{convection_diffusion_control_functional}--\eqref{convection_diffusion_control_constraints} satisfies:
\begin{equation}\label{convection_diffusion_control_optimality_conditions}
\left\{
\begin{array}{rl}

\left.
\begin{array}{rl}
\vspace{1ex}
\dfrac{\partial y}{\partial t} - \epsilon\nabla^2 y +\mathbf{w} \cdot \nabla y = \dfrac{1}{\beta}\: p + f& \quad \mathrm{in} \; \Omega \times (0,t_{f}) \\
\vspace{1ex}
y(\mathbf{x},t) = g(\mathbf{x},t) & \quad \mathrm{on} \; \partial \Omega \times (0,t_{f})\\
\vspace{1ex}
y(\mathbf{x},0) = y_0(\mathbf{x}) & \quad \mathrm{in} \; \Omega
\end{array}
 \right\} & 
\left.
\begin{array}{c}
\vspace{0.5ex}
\mathrm{state}\\
\mathrm{equation}
\end{array}
\right.

\\

\left.
\begin{array}{rl}
\vspace{1ex}
-\dfrac{\partial p}{\partial t} - \epsilon\nabla^2 p - \mathbf{w} \cdot \nabla p \: = \: \widehat{y} - y  & \quad \mathrm{in} \; \Omega \times (0,t_{f}) \\
\vspace{1ex}
p(\mathbf{x},t) = 0 & \quad \mathrm{on} \; \partial \Omega \times (0,t_{f})\\
\vspace{1ex}
p(\mathbf{x},t_{f})= 0 & \quad \mathrm{in} \; \Omega
\end{array}
\right\} & 
\left.
\begin{array}{c}
\vspace{0.5ex}
\mathrm{adjoint}\\
\mathrm{equation}
\end{array}
\right.

\end{array}
\right.
\end{equation}
where we have substituted the gradient equation $\beta u - p = 0$ into the state equation.

Problem \eqref{convection_diffusion_control_optimality_conditions} is a coupled system of (time-dependent) PDEs, consisting of a forward PDE combined with a backward problem for the adjoint. Due to this structure, when trying to find numerical approximation of the solution, an A--stable method is usually applied for discretizing the time derivative, since this will allow the user to choose a time step that is independent of the spatial mesh-size. Further, in order to obtain a consistent system of linear equations, both functions $y$ and $p$ are approximated at the same time points.

We now introduce methods for approximating the time derivative when solving \eqref{convection_diffusion_control_optimality_conditions}, and describe the resulting matrix systems, starting with the widely-used backward Euler method in Section \ref{3_1}, followed by the Crank--Nicolson method in Section \ref{3_2}. For the remainder of the paper, we discretize the interval $(0,t_{f})$ into $n_t$ subintervals of length $\tau=\frac{t_{f}}{n_t}$, and we use the notation $\boldsymbol{y}^n \approx y(\mathbf{x},t_n)$, $\boldsymbol{p}^n \approx p(\mathbf{x},t_n)$ for our approximations for all $\mathbf{x} \in \Omega$, with $t_n=n\tau$.

\subsection{Backward Euler discretization}
\label{3_1}
Many widely-used preconditioned iterative methods for the solution of time-dependent PDE-constrained optimization problems of type \eqref{convection_diffusion_control_functional}--\eqref{convection_diffusion_control_constraints} involve a backward Euler discretization for the time variable \cite{Dolgov_Stoll, Pearson_Stoll, Pearson_Stoll_Wathen, Stoll_Benner_Onwunta_Dolgov, Stoll_Pearson_Maini, Yucel_Stoll_Benner}. Applying this scheme gives the following system of equations\footnote{Note that the matrix $(L^{E})^\top$ arises in the discretization due to the skew-symmetry of $N$, which holds as Assumption \ref{assumption_w} is trivially satisfied for Dirichlet problems.}:
\begin{displaymath}
\left\{
\begin{array}{rl}
\vspace{1ex}
\displaystyle 
\tau M \boldsymbol{y}^n + (L^{E})^\top \boldsymbol{p}^n - M \boldsymbol{p}^{n+1} = \tau M \widehat{\boldsymbol{y}}^n ,
& \quad n=0,...,n_t-1,\\
\vspace{1ex}
\boldsymbol{p}^{n_t}=\mathbf{0},\\
\vspace{1ex}
\boldsymbol{y}^{0}=\boldsymbol{y}_0,\\
\displaystyle - M \boldsymbol{y}^n + L^{E} \boldsymbol{y}^{n+1} - \frac{\tau}{\beta}M \boldsymbol{p}^{n+1}  = \tau \boldsymbol{f}^{n+1},
& \quad n=0,...,n_t-1,\\
\end{array}
\right.
\end{displaymath}
where $L^{E}= \tau L + M$, $\boldsymbol{y}_0$ is the discretization of the initial condition for $y$, and
\begin{equation}\label{discretized_rhs}
\boldsymbol{f}^n=\{\boldsymbol{f}^n_i\}_{i=1}^{n_x}, \qquad \boldsymbol{f}^n_i = \int_{\Omega} f(\mathbf{x},t_n) \phi_i \: \mathrm{d} \Omega \: .
\end{equation}
With this notation, we obtain the following (symmetric) linear system:
\begin{equation}\label{heat_control_discretized}
\left[
\begin{array}{cc}\vspace{1ex}
A^{E} & (B^{E})^\top\\
B^E & -C^E
\end{array}
\right]
\left[
\begin{array}{c}\vspace{1ex}
\bar{\mathbf{y}}\\
\bar{\mathbf{p}}
\end{array}
\right]
= 
\left[
\begin{array}{c}\vspace{1ex}
\textbf{b}_1^{E} \\
\textbf{b}_2^{E}
\end{array}
\right],
\end{equation}
where
\begin{align*}
&A^E = \left[\begin{array}{cccc}
\tau M & & & \\ & \ddots & & \\ & & \tau M & \\ & & & 0 \\
\end{array}\right], \qquad B^E = \left[\begin{array}{cccc}
L^{E} & & & \\ -M & L^{E} & & \\ & \ddots & \ddots & \\ & & -M & L^{E} \\
\end{array}\right], \\
&C^E = \left[\begin{array}{cccc}
0 & & & \\ & \frac{\tau}{\beta} M & & \\ & & \ddots & \\ & & & \frac{\tau}{\beta} M \\
\end{array}\right], \quad \left. \begin{array}{c} \bar{\mathbf{y}} = \left[\begin{array}{cccc} (\boldsymbol{y}^0)^\top & (\boldsymbol{y}^1)^\top & \hdots & (\boldsymbol{y}^{n_t} )^\top \end{array} \right]^\top, \\ \\ \bar{\mathbf{p}} = \left[\begin{array}{cccc} (\boldsymbol{p}^0)^\top & (\boldsymbol{p}^1)^\top & \hdots & (\boldsymbol{p}^{n_t} )^\top \end{array} \right]^\top, \end{array} \right.
\end{align*}
and the right-hand side vectors $\textbf{b}_1^{E}$ and $\textbf{b}_2^{E}$ take into account the initial and final-time conditions on $y$ and $p$, as well as information from the desired state $\widehat{y}$ and the force function $f$. Note that we have symmetrized the system by replacing the discretized initial and final-time conditions with
\begin{displaymath}
(L^{E})^\top \boldsymbol{p}^{n_t}=\mathbf{0}, \qquad L^{E} \boldsymbol{y}^{0}= L^{E} \boldsymbol{y}_0.
\end{displaymath}
The structure of the previous system is very convenient from the point of view of numerical linear algebra because it facilitates effective preconditioning.

\subsection{Crank--Nicolson discretization and symmetrization of the system}
\label{3_2}
Despite the convenient structure of the matrix system arising from the backward Euler method, the essential drawback is that the method is only first-order accurate in time. For instance, if a second-order accurate method is used for the spatial discretization, the numerical error is of order $\mathcal{O}(\tau) + \mathcal{O}(h^2)$, given sufficient smoothness properties of the solution. Therefore, it is reasonable to choose $\tau = \mathcal{O}(h^2)$, resulting in a matrix system of huge dimension, and hence an extremely high CPU time is required for its solution. In recent years, a significant effort has been invested in improving the accuracy of the discretized solution of time dependent PDE-constrained optimization problems involving the backward Euler method: see \cite{Guttel_Pearson} for an application of deferred correction to time-dependent PDE-constrained optimization problems, \cite{Dolgov_Stoll, Stoll, Stoll_Benner_Onwunta_Dolgov, Stoll_Breiten} for low-rank tensor approaches to speed up the convergence of backward Euler, \cite{Hinze_Koster_Turek} for a multigrid technique applied to optimal flow control problems, and \cite{Zeng_Zhang} for a preconditioned iterative solver for problem of the type \eqref{heat_control_functional}--\eqref{heat_control_constraints} that uses a reduction in dimensionality of the arising system.

We also highlight recent works in which higher-order time discretizations are considered. For example, in \cite{Gotschel_Minion} the authors derive a parallel-in-time algorithm coupled with a gradient-based method, using a Lobatto IIIA scheme in the time variable, and in \cite{Dolgov_Stoll} the authors note that their low-rank method can also be adopted to the Crank--Nicolson approach. Other valuable approaches for addressing time-dependent PDE-constrained optimization problems include multiple shooting methods \cite{Hein05}, parareal schemes \cite{MaTu02,MSS10}, and ideas from instantaneous control \cite{CHK99,Hinz00}. However, to the best of our knowledge, the crucial question of preconditioning PDE-constrained optimization methods by exploiting the precise matrix structures arising from higher-order discretization schemes in time has not been resolved, perhaps due to the increased complexity in the structure of the matrix systems, an issue pointed out in \cite{Hinze_Koster_Turek} for instance.

\begin{remark}
We note here that both backward Euler and Crank--Nicolson are unconditionally A--stable, that is, no restriction on $\tau$ and $h$ is required. There are often other stability properties to consider when selecting an appropriate time-stepping scheme: for instance Crank--Nicolson, as opposed to backward Euler, is not L--stable, which is an important consideration when applying long-range integrations. In this work, we do not address specific questions about stability properties of the two methods, but rather whether it is possible in principal to take advantage of a higher-order discretization scheme in the time variable for increasing the rate of convergence of solvers for the optimal control of time-dependent PDEs. It is likely that our proposed method could be extended to time-stepping methods that achieve even faster convergence than Crank--Nicolson, including L--stable methods. We highlight the reference \cite{AxBlKo}, which considers the question of preconditioning matrix systems arising from L--stable methods for forward evolutionary PDEs.
\end{remark}

\begin{remark}
We note that, in order for Crank--Nicolson to achieve a second-order convergence rate, we require the state $y$ and the adjoint variable $p$ in \eqref{convection_diffusion_control_optimality_conditions} to have sufficient smoothness. Therefore, in the following analysis, we allow the assumptions of regularity and the results in \cite[Sec.\,3]{Apel_Flaig} to hold for the problem considered.
\end{remark}

Considering again \eqref{convection_diffusion_control_optimality_conditions}, we now discretize the time derivative using the Crank--Nicolson method. Denoting
\begin{equation*}
L^+ = \frac{\tau}{2} L + M, \quad  L^- = \frac{\tau}{2} L - M, \quad \bar{M}= \frac{\tau}{2}M, \quad \bar{M}_\beta=\frac{\tau}{2 \beta}M,
\end{equation*}
we have that the numerical solution of \eqref{convection_diffusion_control_optimality_conditions} satisfies
\begin{displaymath}
\left\{
\begin{array}{rl}
\displaystyle 
\bar{M} \left( \boldsymbol{y}^n + \boldsymbol{y}^{n+1} \right) + (L^+)^\top \boldsymbol{p}^n + (L^-)^\top \boldsymbol{p}^{n+1} = \bar{M} \left(\widehat{\boldsymbol{y}}^n + \widehat{\boldsymbol{y}}^{n+1}\right),
& n=0,...,n_t-1,\\
\\
\displaystyle  L^- \boldsymbol{y}^n + L^+ \boldsymbol{y}^{n+1} - \bar{M}_\beta \left(\boldsymbol{p}^n + \boldsymbol{p}^{n+1} \right) =\frac{\tau}{2} \left(\boldsymbol{f}^n + \boldsymbol{f}^{n+1}\right),
&  n=0,...,n_t-1,\\
\end{array}
\right.
\end{displaymath}
with $M\boldsymbol{p}^{n_t} = \mathbf{0}$ and $M\boldsymbol{y}^0 = M\boldsymbol{y}_0$ an appropriate discretization of the final and initial conditions on $p$ and $y$, and $\boldsymbol{f}^n$ defined as in \eqref{discretized_rhs}. In matrix form, we write
\begin{equation}\label{original_system}
\left.
\begin{array}{c}

\left[
\begin{array}{cc}
\vspace{1ex}
\bar{\Lambda}_{11} & \bar{\Lambda}_{12}\\

\bar{\Lambda}_{21} & \bar{\Lambda}_{22}\\
\end{array}

\right]

\left[
\begin{array}{c}
\vspace{1ex}
\bar{\mathbf{y}} \\
\bar{\mathbf{p}}
\end{array}
\right]
=
\left[
\begin{array}{c}
\vspace{1ex}
\bar{\mathbf{b}}_1 \\
\bar{\mathbf{b}}_2
\end{array}
\right],
\end{array}
\right.
\end{equation}
where the vectors $\bar{\mathbf{y}}$ and $\bar{\mathbf{p}}$ are the numerical solution for the state and adjoint variables, and as before the right-hand side accounts for the initial and final-time conditions on $y$ and $p$, as well as the desired state $\widehat{y}$ and force function $f$. The matrices $\bar{\Lambda}_{ij}$, $i,j=1,2$, are given by
\begin{align*}
&\bar{\Lambda}_{11} = 
\left[
\begin{array}{cccc}
\bar{M} & \bar{M} &   \\
 & \ddots & \ddots & \\
  & & \bar{M} & \bar{M}\\
 &  &  & 0
\end{array}
\right], \quad \bar{\Lambda}_{12} = 
\left[
\begin{array}{cccc}
(L^+)^\top & (L^-)^\top &  \\
 & \ddots & \ddots & \\
  &  & (L^+)^\top & (L^-)^\top \\
  &  &  & M
\end{array}
\right], \\
&\bar{\Lambda}_{21} = 
\left[
\begin{array}{cccc}
M & \\
L^- & L^+ & \\
 & \ddots & \ddots & \\
  &  & L^- & L^+
\end{array}
\right], \quad \bar{\Lambda}_{22} = -
\left[
\begin{array}{cccc}
0 & \\
\bar{M}_\beta & \bar{M}_\beta &  \\
 & \ddots & \ddots & \\
 &  & \bar{M}_\beta & \bar{M}_\beta \\
\end{array}
\right].
\end{align*}
We immediately realize that the matrices $\bar{\Lambda}_{11}$ and $\bar{\Lambda}_{22}$ are not symmetric, and therefore no iterative method for symmetric matrices, such as MINRES, may be applied to \eqref{original_system}. We therefore wish to apply a transformation to \eqref{original_system} to convert it to a symmetric system. We partition the matrix in \eqref{original_system} as

\begin{displaymath}
\left.
\begin{array}{c}

\left[
\begin{array}{c|cccc|cccc|c}
\bar{M} & \bar{M} &  & & & (L^+)^\top & (L^-)^\top & & &0 \\
0 & \bar{M} & \bar{M} & & & & (L^+)^\top & (L^-)^\top & \\
 & & \ddots & \ddots & & & & \ddots & \ddots &\\
  & & & \bar{M} & \bar{M}   & & & & (L^+)^\top & (L^-)^\top \\
  \hline
 & & & &  0 & & &  &  & M\\
\hline
M & & & & & 0 & & & \\
\hline
L^- & L^+ & & & & -\bar{M}_\beta & -\bar{M}_\beta &  & & 0 \\
& L^- & L^+ & & & & -\bar{M}_\beta & -\bar{M}_\beta &  \\
 & & \ddots & \ddots &  & & & \ddots & \ddots & \\
  &  &&  L^- & L^+ &  & &&-\bar{M}_\beta & -\bar{M}_\beta \\
\end{array}

\right].
\end{array}
\right.
\end{displaymath}
Then if we eliminate the initial and final-time conditions on $y$ and $p$, we can rewrite
\begin{displaymath}
\left.
\begin{array}{c}

\left[
\begin{array}{cc}
\vspace{1ex}
\Lambda_{11} & \Lambda_{12}\\

\Lambda_{21} & -\Lambda_{22}\\
\end{array}

\right]

\left[
\begin{array}{c}
\vspace{1ex}
\mathbf{y} \\
\mathbf{p}
\end{array}
\right]
=
\left[
\begin{array}{c}
\vspace{1ex}
\mathbf{b}_1 \\
\mathbf{b}_2
\end{array}
\right],
\end{array}
\right.
\end{displaymath}
with the vectors $\mathbf{y}$, $\mathbf{p}$, $\mathbf{b}_1$, $\mathbf{b}_2$ modified accordingly, and the matrices $\Lambda_{ij}$, $i,j=1,2$, given by
\begin{align}
\nonumber &\Lambda_{11} = 
\left[
\begin{array}{cccc}
\bar{M}\\
\bar{M} & \bar{M} &   \\
 & \ddots & \ddots & \\
  & & \bar{M} & \bar{M}\\
\end{array}
\right] \quad \Lambda_{12} = 
\left[
\begin{array}{cccccc}
(L^+)^\top & (L^-)^\top &  \\
 & \ddots & \ddots &  \\
  &  & (L^+)^\top & (L^-)^\top \\
  &  &  & (L^+)^\top
\end{array}
\right], \\
\label{Lambda_21} &\Lambda_{21}=
\underbrace{\left[
\begin{array}{cccccc}
L^+ & \\
L^- & L^+ & \\
 & \ddots & \ddots & \\
  &  & L^- & L^+
\end{array}
\right]}_{ = \Lambda_{12} ^\top}, \quad \Lambda_{22} = 
\left[
\begin{array}{cccc}
\bar{M}_\beta & \bar{M}_\beta &  \\
 & \ddots & \ddots & \\
 &  & \bar{M}_\beta & \bar{M}_\beta \\
 & & & \bar{M}_\beta
\end{array}
\right].
\end{align}

In order to symmetrize the system, we now apply the following linear transformation:
\begin{equation}\label{linear_transf_T}
T=
\left[
\begin{array}{cc}\vspace{1ex}
T_1 & 0  \\
0 & T_2
\end{array}
\right],
\quad T_1=
\left[
\begin{array}{cccc}
I & I &  \\
 & \ddots & \ddots  \\
 & & I & I\\
 &  &  & I
\end{array}
\right],
\quad T_2=
\underbrace{\left[
\begin{array}{cccc}
I &  \\
I & I &  \\
 & \ddots & \ddots  \\
 &   & I & I
\end{array}
\right]}_{= T_1^{\top}},
\end{equation}
where $T_1,T_2 \in \mathbb{R}^{(n_t n_x) \times (n_t n_x)}$, and $I\in \mathbb{R}^{n_x \times n_x}$ denotes the identity matrix. Then,
\begin{equation}\label{modified_system}
T\left[
\begin{array}{cc}\vspace{1ex}
\Lambda_{11} & \Lambda_{12}\\
\Lambda_{21} & \Lambda_{22}
\end{array}
\right]
\left[
\begin{array}{c}\vspace{1ex}
\mathbf{y}\\
\mathbf{p}
\end{array}
\right]
=
\underbrace{\left[
\begin{array}{cc}\vspace{1ex}
A & B^\top\\
B & -C
\end{array}
\right]}_{\mathcal{A}}
\left[
\begin{array}{c}\vspace{1ex}
\mathbf{y}\\
\mathbf{p}
\end{array}
\right]
=
T
\left[
\begin{array}{c}\vspace{1ex}
\mathbf{b}_1 \\
\mathbf{b}_2
\end{array}
\right],
\end{equation}
where
\begin{equation*}
A = T_1 \Lambda_{11} = 
\left[
\begin{array}{cccc}
2\bar{M} & \bar{M} & & \\
\bar{M} & \ddots & \ddots &   \\
  & \ddots & 2\bar{M} & \bar{M}\\
 &  & \bar{M} & \bar{M}
\end{array}
\right],
\quad \! C = T_2 \Lambda_{22} = 
\left[
\begin{array}{cccc}
\bar{M}_\beta & \bar{M}_\beta & & \\
\bar{M}_\beta & 2\bar{M}_\beta & \ddots & \\
 & \ddots & \ddots & \bar{M}_\beta \\
 & & \bar{M}_\beta & 2\bar{M}_\beta\\
\end{array}
\right].
\end{equation*}
Furthermore, it holds that
\begin{equation*}
B=T_2 \Lambda_{21} = 
\left[
\begin{array}{cccccc}
L^+ & \\
\tau L & L^+ &  \\
L^- & \tau L & L^+ \\
 & \ddots & \ddots & \ddots & \\
   & & L^- & \tau L & L^+
\end{array}\right]=(T_1 \Lambda_{12})^\top.
\end{equation*}
We have thus transformed \eqref{original_system} to a symmetric system using matrices $T_1$, $T_2$ which are easy and cheap to apply and invert. We also highlight the further properties that
\begin{equation}\label{A_C_transform}
A=T_1 A_D T_1^\top, \qquad
C=T_2 C_D T_2^\top,
\end{equation}
where
\begin{equation}\label{bar_AC}
A_D
=  
\left[
\begin{array}{ccccc}
\bar{M} &   \\
 & \ddots &   \\
 &   & \bar{M}\\
\end{array}
\right],
\qquad
C_D
=  
\left[
\begin{array}{ccccc}
\bar{M}_\beta &   \\
 & \ddots &   \\
 &   & \bar{M}_\beta\\
\end{array}
\right],
\end{equation}
so we may work with the matrices $A$ and $C$ cheaply, using $T_1$, $T_2$, $A_D$, $C_D$. Further, since both $A_D$ and $C_D$ are symmetric positive definite, the same holds for $A$ and $C$.

So, in order to find an approximate solution to \eqref{convection_diffusion_control_optimality_conditions}, we may now consider the \emph{saddle-point system} \eqref{modified_system}, to which we can apply a preconditioned Krylov subspace method for symmetric indefinite matrices, such as MINRES.

\section{Preconditioning Approach}
\label{4}
In this section we describe an optimal preconditoner for the system \eqref{modified_system}, by making use of saddle-point theory as well as suitable approximations for the blocks of the matrix $\mathcal{A}$.

\subsection{Saddle-point systems}
\label{4_1}
We know that if we wish to solve an invertible system of the form \eqref{modified_system}, with invertible $A$, a good candidate for a preconditioner is the block diagonal matrix:
\begin{equation}\label{optimal_prec}
\mathcal{P}=
\left[
\begin{array}{cc}
A & 0 \\
0 & S
\end{array}
\right],
\end{equation}
where $S$ denote the (negative) \emph{Schur complement}  $S=C+BA^{-1}B^\top$. Indeed, if $C=0$ and $S$ is also invertible, the preconditioned matrix $\mathcal{P}^{-1}\mathcal{A}$ will have only 3 distinct eigenvalues \cite{Kuznetsov:1995:EIS, Murphy:1999:NPI:359189.359190}:
\begin{displaymath}
\lambda (\mathcal{P}^{-1}\mathcal{A}) = \left\{1, \dfrac{1\pm \sqrt{5}}{2} \right\},
\end{displaymath}
where we denote the set of eigenvalues of a matrix by $\lambda(\cdot)$. If instead $A$ and $C$ are symmetric positive definite, it holds that (see \cite{AxNe06,SiWa94} and \cite[Theorem 4]{PearsonPhD}):
\begin{displaymath}
\lambda (\mathcal{P}^{-1}\mathcal{A}) \in \left[-1,\frac{1}{2}(1-\sqrt{5})\right] \cup \left[1,\frac{1}{2}(1+\sqrt{5})\right].
\end{displaymath}

In these cases, a preconditioned Krylov subspace method for symmetric indefinite matrices with preconditioner \eqref{optimal_prec} would converge in few iterations (at most 3 if $C=0$). We further observe that, if $A$ and $C$ are symmetric positive definite, the preconditioner $\mathcal{P}$ also has that property, so a natural choice of iterative solver is the MINRES algorithm. Of course, as in \eqref{optimal_prec} the preconditioner $\mathcal{P}$ is not practical, as the exact inverses of $A$ and $S$ will need to be applied in each case. Applying $S^{-1}$ is particularly problematic as even when $A$ and $C$ are sparse, $S$ is generally dense. For this reason, we wish to find a suitable approximation $\widehat{\mathcal{P}}$ of $\mathcal{P}$, with
\begin{displaymath}
\widehat{\mathcal{P}} = 
\left[
\begin{array}{cc}
\widehat{A} & 0 \\
0 & \widehat{S}
\end{array}
\right],
\end{displaymath}
or, more precisely, a cheap application of the effect of $\widehat{\mathcal{P}}^{-1}$ on a generic vector. In the following section we commence by finding a good approximation of the inverse of the matrix $A$, then in Section \ref{4_3} we describe the approach used for approximating the Schur complement $S$.

Since we will benchmark our new preconditioning strategy against the preconditioner derived in \cite{Pearson_Stoll_Wathen} for heat control (i.e., $L=K$) after applying the optimize-then-discretize approach with backward Euler time-stepping, we briefly describe the preconditioner for the system \eqref{heat_control_discretized}, itself also a saddle-point system. In this case the $(1,1)$-block is not invertible; a good preconditioner is found to be
\begin{displaymath}
\widehat{\mathcal{P}}^E = 
\left[
\begin{array}{cc}
\widehat{A}^E & 0 \\
0 & \widehat{S}^E
\end{array}
\right],
\end{displaymath}
where
\begin{equation}\label{prec_Euler}
\widehat{A}^E=\left[
\begin{array}{cccc}
\tau  M & &\\
  & \ddots & \\
 &  &\tau M & \\
 &   &  & \xi \tau M
\end{array}
\right], \qquad \widehat{S}^E = (B^E + \widehat{M}^E) (\widehat{A}^E)^{-1}  (B^E + \widehat{M}^E)^\top,
\end{equation}
with
\begin{equation}\label{hat_M_Euler}
\widehat{M}^E = \dfrac{\tau}{\sqrt{\beta}} 
\left[
\begin{array}{cccccc}
M & \\
 & \ddots & \\
 &  & M\\
 & & & \sqrt{\xi} M
\end{array}
\right].
\end{equation}
Here, $0<\xi\ll 1$ is chosen such that the (invertible) matrix $\widehat{A}^E$ is `close enough' to the matrix $A^E$ in some sense. The $(1,1)$-block is approximated by a `perturbed' matrix $\widehat{A}^E$, whose inverse is then applied within the matrix $\widehat{S}^E$. A term of the form \eqref{hat_M_Euler} is justified in Section \ref{4_3}.

\subsection{Approximation of the \boldmath{$(1,1)$}-block}
\label{4_2}
We now focus on devising a preconditioner for the matrix system \eqref{modified_system}, arising from a Crank--Nicolson discretization, starting with a cheap and effective approximation of the matrix $A$. We recall from \eqref{A_C_transform} that $A$ can be written as $A=T_1 A_D T_1^\top$, with $A_D$ defined as in \eqref{bar_AC}. We observe that $A_D$ is a block diagonal matrix with each diagonal block a multiple of $M$. Therefore, in order to obtain a cheap approximation for $A$, we require a suitable way to approximate the inverse of a mass matrix. As discussed in \cite{GolubVargaI,GolubVargaII,Wathen_Rees}, the Chebyshev semi-iterative method represents a cheap and effective way to do this. From the discussion above, a good approximation of $A$ is therefore given by $\widehat{A}=T_1 \widehat{A}_D T_1^\top$, with
\begin{displaymath}
\widehat{A}_D = \dfrac{\tau}{2}
\left[
\begin{array}{ccccc}
M_c &   \\
 & \ddots &   \\
 &   & M_c\\
\end{array}
\right],
\end{displaymath}
where $M_c$ denotes a fixed number of steps (20, in our tests) of Chebyshev semi-iteration applied to $M$. Further, since $T_1$ and $T_2 = T_1^\top$ are (block) upper- and lower-triangular matrices respectively (which makes them very easy to invert), we use backward and forward (block-)substitution in order to apply their inverses.

It is trivial to prove the following, on the effectiveness of our approximation of $A$:
\begin{lemma}
The minimum (maximum) eigenvalue of $\widehat{A} ^{-1} A$ is bounded below (above) by the minimum (maximum) eigenvalue of $M_c^{-1} M$.
\end{lemma}

\subsection{Approximation of Schur complement}
\label{4_3}
We now find a suitable approximation for the Schur complement $S$ of \eqref{modified_system}. In the forthcoming theory, we assume that Assumption \ref{assumption_w} holds, and later discuss the case where this is relaxed. Recalling how the matrices $A$, $B$, and $C$ can be written as linear transformations involving the matrices $T_1$ and $T_2$, we first rewrite $S$ in the following way:
\begin{equation}\label{schur_complement}
S = T_2 C_D T_2^\top + T_2 \Lambda_{21} (T_1 A_D T_1^\top)^{-1} \Lambda_{21}^\top T_2^\top=T_2\left[ C_D + B_D A_D^{-1} B_D^\top \right] T_2^\top,
\end{equation}
where we set
\begin{equation}\label{B_D}
B_D = \Lambda_{21}T_2^{-1}.
\end{equation}
It is hence clear that if we find a symmetric positive definite approximation $\widetilde{S}_{\mathrm{int}}$ of
\begin{equation}\label{S_int_def}
S_{\mathrm{int}}=C_D + B_D A_D^{-1} B_D^\top,
\end{equation}
then $\widehat{S}:=T_2 \widetilde{S}_{\mathrm{int}} T_2^\top$ is a symmetric positive definite approximation of $S$. We may therefore use the (generalized) Rayleigh quotient in order to find an upper and lower bounds for the eigenvalues of the matrix $\widehat{S}^{-1}S$ as follows. Let $\lambda$ be an eigenvalue of $\widehat{S}^{-1}S$ with $\mathbf{x}$ the corresponding eigenvector; then
\begin{equation}\label{Rayleigh_quotient}
\widehat{S}^{-1}S \mathbf{x} = \lambda \mathbf{x} \; \Rightarrow \; S \mathbf{x} = \lambda  \widehat{S}\mathbf{x} \; \Rightarrow \;
\lambda = \dfrac{\mathbf{x}^{\top}S \mathbf{x}}{\mathbf{x}^{\top}\widehat{S} \mathbf{x}} = \dfrac{\widetilde{\mathbf{x}}^{\top}S_{\mathrm{int}} \widetilde{\mathbf{x}}}{\widetilde{\mathbf{x}}^{\top}\widetilde{S}_{\mathrm{int}} \widetilde{\mathbf{x}}},
\end{equation}
where we set $\widetilde{\mathbf{x}}=T_2^\top \mathbf{x}\neq \mathbf{0}$ and use \eqref{schur_complement} together with the definition of $\widehat{S}$. Thus upper and lower bounds for the eigenvalues of $\widehat{S}^{-1}S$ are given by the maximum and the minimum eigenvalues of $\widetilde{S}_{\mathrm{int}}^{-1}S_{\mathrm{int}}$, respectively.

From \eqref{schur_complement}, and due to the convenient structure of the matrices $A_D$ and $C_D$ in \eqref{bar_AC}, we can devise an approximation $\widetilde{S}_{\mathrm{int}}$ of $ S_{\mathrm{int}}$ using the `matching strategy' discussed in \cite{Pearson_Stoll_Wathen,Pearson2018} as follows. We seek an approximation:
\begin{equation}\label{tilde_S_int_def}
\widetilde{S}_{\mathrm{int}} = (\Lambda_{21} + \widehat{M} )A^{-1} (\Lambda_{21} + \widehat{M})^\top \approx S_{\mathrm{int}} 
\end{equation}
such that $\widetilde{S}_{\mathrm{int}}$ `captures' both terms of $S_{\mathrm{int}}$, namely $C_D$ and $B_D A_D^{-1} B_D^\top$ (noting that $\Lambda_{21}A^{-1}\Lambda_{21}^\top=B_D A_D^{-1} B_D^\top$). We therefore wish that
\begin{equation}\label{appr_C_D}
\widehat{M} A^{-1} \widehat{M}^\top = \left[ \widehat{M} (T_1^\top)^{-1} \right] A_D^{-1} \left[T_1^{-1} \widehat{M}^\top\right]  = C_D.
\end{equation}
Using the definitions of $C_D$ and $A_D$ from \eqref{bar_AC}, we obtain that for this to hold the matrix $\widehat{M} (T_1^\top)^{-1}$ is given by
\begin{equation}\label{M_D}
M_D := \widehat{M} (T_1^\top)^{-1} = \frac{\tau}{2\sqrt{\beta}}
\left[
\begin{array}{cccccc}
M & \\
   & \ddots & \\
  & & M 
\end{array}
\right],
\end{equation}
and therefore
\begin{equation}\label{Mhat}
\widehat{M} = \dfrac{\tau}{2\sqrt{\beta}} 
\left[
\begin{array}{cccccc}
M & \\
  M & M & \\
   & \ddots & \ddots & \\
  & & M & M
\end{array}
\right].
\end{equation}

Finally, our approximation of $S$ is given by
\begin{equation}\label{hat_S_def}
\widehat{S} = T_2 (\Lambda_{21} + \widehat{M} )A^{-1} (\Lambda_{21} + \widehat{M})^\top T_2^\top = (\Lambda_{21} + \widehat{M} )A_D^{-1} (\Lambda_{21} + \widehat{M}),
\end{equation}
with $\widehat{M}$ as defined in \eqref{Mhat}, and the two expressions are equivalent since $T_2$ commutes with both $\Lambda_{21}$ and $\widehat{M}$. To understand the effectiveness of this approximation, we recall \eqref{Rayleigh_quotient}, telling us that we only need to study the spectrum of the matrix $\widetilde{S}_{\mathrm{int}}^{-1}S_{\mathrm{int}}$. We next rewrite $\widetilde{S}_{\mathrm{int}}$ as follows:
\begin{align}
\nonumber \widetilde{S}_{\mathrm{int}} ={}& (\Lambda_{21} + \widehat{M} )A^{-1} (\Lambda_{21} + \widehat{M})^\top \\
\nonumber ={}& \widehat{M} A^{-1} \widehat{M}^{\;\top} + \Lambda_{21} A^{-1} \Lambda_{21}^\top + 
 \widehat{M} A^{-1}\Lambda_{21}^\top + \Lambda_{21} A^{-1} \widehat{M}^{\;\top  }\\
\nonumber ={}& C_D + B_D A_D^{-1} B_D^\top + 
 M_D A_D^{-1}B_D^\top + B_D A_D^{-1} M_D^\top \\
\label{tilde_S} ={}& S_{\mathrm{int}} + 
 M_D A_D^{-1}B_D^\top + B_D A_D^{-1} M_D^\top ,
\end{align}
where we have used \eqref{appr_C_D}, \eqref{A_C_transform}, \eqref{B_D}, and \eqref{M_D} in turn.

Since $S_{\mathrm{int}}$ and $\widetilde{S}_{\mathrm{int}}$ are symmetric positive definite, we again consider the generalized Rayleigh quotient:
\begin{equation}\label{Rayleigh_quotient_S_int}
R : = \dfrac{\mathbf{x}^{\top}S_{\mathrm{int}} \mathbf{x}}{\mathbf{x}^{\top}\widetilde{S}_{\mathrm{int}} \mathbf{x}} = \dfrac{\mathbf{a}^\top \mathbf{a} + \mathbf{b}^\top \mathbf{b}}{\mathbf{a}^\top \mathbf{a} + \mathbf{b}^\top \mathbf{b} + \mathbf{a}^\top \mathbf{b} + \mathbf{b}^\top \mathbf{a}},
\end{equation}
where $\mathbf{a}= (B_D A_D^{-1/2})^\top \mathbf{x}$ and $\mathbf{b}= (C_D^{1/2})^\top \mathbf{x}$, noting from \eqref{bar_AC} and \eqref{M_D} that $C_D ^{1/2} = M_D A_D^{-1/2}$. As $\mathbf{b}^\top \mathbf{b} > 0$, we have
\begin{equation*}
R \geq \dfrac{1}{2} ~~ \Leftrightarrow ~~ \mathbf{a}^\top \mathbf{a} + \mathbf{b}^\top \mathbf{b} \geq \dfrac{1}{2}\left( \mathbf{a}^\top \mathbf{a} + \mathbf{b}^\top \mathbf{b} + \mathbf{a}^\top \mathbf{b} + \mathbf{b}^\top \mathbf{a} \right) ~~ \Leftrightarrow ~~ \frac{1}{2} \left(\mathbf{a} - \mathbf{b} \right)^\top \left(\mathbf{a} - \mathbf{b} \right) \geq 0,
\end{equation*}
which is clearly satisfied since $\left(\mathbf{a} - \mathbf{b} \right)^\top\left(\mathbf{a} - \mathbf{b} \right) = \| \mathbf{a} - \mathbf{b} \|^2$.

In order to find an upper bound for the Rayleigh quotient \eqref{Rayleigh_quotient_S_int}, we return to \eqref{tilde_S}. Noting that
\begin{displaymath}
M_D A_D ^{-1} = \dfrac{1}{\sqrt{\beta}} \left[
\begin{array}{lll}
I \\
& \ddots \\
 & & I
\end{array}
\right],
\end{displaymath}
we can rewrite
\begin{displaymath}
\widetilde{S}_{\mathrm{int}} = S_\mathrm{int} + \dfrac{1}{\sqrt{\beta}} \left(
B_D^\top + B_D \right).
\end{displaymath}
Following the reasoning in \cite{Pearson_Stoll_Wathen}, we can prove that $R\leq 1$; from \eqref{Rayleigh_quotient_S_int}, this holds if
\begin{displaymath}
\begin{array}{ll}
\vspace{1ex}
\mathbf{a}^\top \mathbf{b} + \mathbf{b}^\top \mathbf{a} \hspace{-0.75em} & = \mathbf{x}^\top \left( M_D A_D^{-1} B_D^\top  + B_D A_D^{-1} M_D^\top \right)\mathbf{x}  \\
\vspace{1ex}
 & = \dfrac{1}{\sqrt{\beta}} \, \mathbf{x}^\top \left(  B_D^\top  + B_D \right)\mathbf{x} \\
 & = \dfrac{1}{\sqrt{\beta}} \, \mathbf{z}^\top \left(  \Lambda_{21}^\top T_2  + T_2^\top \Lambda_{21} \right)\mathbf{z}  \geq 0 \:,
\end{array}
\end{displaymath}
where the last line uses \eqref{B_D} and sets $\mathbf{z} = T_2^{-1} \mathbf{x}$. Therefore, we wish to show that the matrix $\mathcal{B} =\Lambda_{21}^\top T_2  + T_2^\top \Lambda_{21}$ is positive semi-definite. We easily obtain that
\begin{displaymath}
\mathcal{B} = 
\dfrac{\tau}{2}\: \underbrace{
\left[
\begin{array}{cccc}
2 \widetilde{L} & \widetilde{L} & & \\
\widetilde{L} & \ddots & \ddots & \\
 & \ddots & 2 \widetilde{L} & \widetilde{L} \\
 & & \widetilde{L} & \widetilde{L}
\end{array}
\right]
}_{=:\widetilde{\mathcal{L}}}
 +
 \underbrace{
\left[
\begin{array}{cccc}
0 &  \\
 & \ddots &  \\
 & & 0 \\
 & & & 2 M
\end{array}
\right]
}_{=:\widetilde{\mathcal{M}}}
,
\end{displaymath}
with $\widetilde{L}=L + L^\top = 2 (\epsilon K + W)$ since both $K$ and $W$ are symmetric, and $N$ is skew-symmetric due to Assumption \ref{assumption_w}. Furthermore, since $K$ is positive definite in this case, with $W$ positive semi-definite, $\widetilde{L}$ is also positive definite.
Moreover, it is clear that
\begin{displaymath}
\mathbf{z} ^\top \widetilde{\mathcal{L}} \: \mathbf{z} \geq 0 \: \wedge \: \mathbf{z} ^\top \widetilde{\mathcal{M}} \: \mathbf{z} \geq 0  \quad \Rightarrow \quad \mathbf{z} ^\top \mathcal{B}\: \mathbf{z} \geq 0.
\end{displaymath}
The second inequality above is clear, since $\widetilde{\mathcal{M}}$ is a block diagonal matrix with $n_t-1$ zero-blocks and the last block consisting of a positive definite matrix. To prove the inequality $\mathbf{z} ^\top \widetilde{\mathcal{L}} \: \mathbf{z} \geq 0$, we first observe that $\widetilde{\mathcal{L}} = \mathcal{T} \otimes \widetilde{L}$, with
\begin{displaymath}
\mathcal{T} = 
\left[
\begin{array}{cccc}
2 & 1 & & \\
1 & \ddots & \ddots & \\
 & \ddots & 2 & 1 \\
 & & 1 & 1
\end{array}
\right] =
\left[
\begin{array}{cccc}
1 & 1 & & \\
 & \ddots & \ddots & \\
 & & 1 & 1\\
 & & & 1
\end{array}
\right]
\left[
\begin{array}{cccc}
1 & & & \\
1 & 1 & & \\
 & \ddots & \ddots & \\
 & & 1 & 1
\end{array}
\right]
=\mathcal{T}_1^\top \mathcal{T}_1 \:
,
\end{displaymath}
and therefore $\mathcal{T}$ is symmetric positive definite (since $\mathcal{T}_1$ has full rank). Next, we state and use \cite[Theorem 13.12]{Laub} to prove that $\widetilde{\mathcal{L}}$ is positive definite:
\begin{theorem}\label{Kronecker}
Let $X\in \mathbb{R}^{m \times m}$ have eigenvalues $\lambda_i$, $i =1 ,\ldots, m$, and let $Y\in \mathbb{R}^{r \times r}$ have eigenvalues $\mu_j$, $j = 1,\ldots, r$. Then, the $mr$ eigenvalues of $X\otimes Y$ are
\begin{displaymath}
\lambda_1\mu_1, \ldots, \lambda_1\mu_r,\lambda_2\mu_1,\ldots,\lambda_2\mu_r,\ldots, \lambda_m\mu_r.
\end{displaymath}
Moreover, if $\mathbf{x}_1,\ldots, \mathbf{x}_l$ are linearly independent right eigenvectors of $X$ corresponding to $\lambda_1,\ldots, \lambda_l$, $l\leq m$, and $\mathbf{z}_1,\ldots, \mathbf{z}_q$ are linearly independent right eigenvectors of $Y$ corresponding to $\mu_1,\ldots, \mu_q$, $q\leq r$, then $\mathbf{x}_i \otimes \mathbf{z}_j \in \mathbb{R}^{mr}$ are linearly independent right eigenvectors of $X \otimes Y$ corresponding to $\lambda_i\mu_j$, $i=1\ldots, l$, $j=1,\ldots,q$.
\end{theorem}

Since both $\mathcal{T}$ and $\widetilde{L}$ are symmetric positive definite, applying Theorem \ref{Kronecker} with $X=\mathcal{T}$ and $Y=\widetilde{L}$ tells us that the matrix $\widetilde{\mathcal{L}}= \mathcal{T} \otimes \widetilde{L}$ is symmetric positive definite. We therefore infer that $\mathcal{B}$ is symmetric positive definite, and hence that
\begin{displaymath}
\mathbf{a}^\top \mathbf{b} + \mathbf{b}^\top \mathbf{a} \geq 0,
\end{displaymath}
with $\mathbf{a}$ and $\mathbf{b}$ as defined above. Finally, the last inequality guarantees that the Rayleigh quotient $R$ in \eqref{Rayleigh_quotient_S_int} satisfies $R\leq 1$.

We have hence proved the following result:
\begin{theorem}\label{spectrum_S_int}
Let $S_{\mathrm{int}}$ and $\widetilde{S}_{\mathrm{int}}$ be defined as in \eqref{S_int_def} and \eqref{tilde_S_int_def} respectively, with the matrices $A_D$, $B_D$, $C_D$, $\Lambda_{21}$, $A$, $\widehat{M}$ defined as in \eqref{bar_AC}, \eqref{B_D}, \eqref{Lambda_21}, \eqref{A_C_transform}, and \eqref{Mhat}. Then, given Assumption \ref{assumption_w}:
\begin{displaymath}
\lambda(\widetilde{S}_{\mathrm{int}} ^{-1} S_{\mathrm{int}}) \in \left[ \dfrac{1}{2}, 1\right].
\end{displaymath}
\end{theorem}

\begin{figure}[htb]
\centering
\caption{Eigenvalues of $\widetilde{S}_{\mathrm{int}}^{-1} S_{\mathrm{int}}$ for $\beta=10^{-i}$, $i=2,3,4,5$, with $\epsilon=\frac{1}{100}$, $\mathbf{w}=[ 2x_2(1-x_1^2), ~-2x_1(1-x_2^2) ]^\top$ (where $\mathbf{x}=[x_1 ,~ x_2 ]^\top$), on an evenly spaced space-time grid $(-1,1)^2 \times (0,2)$ with $\tau=h=\frac{1}{8}$.}
\label{fig:eigs}
\includegraphics[width=0.49\linewidth]{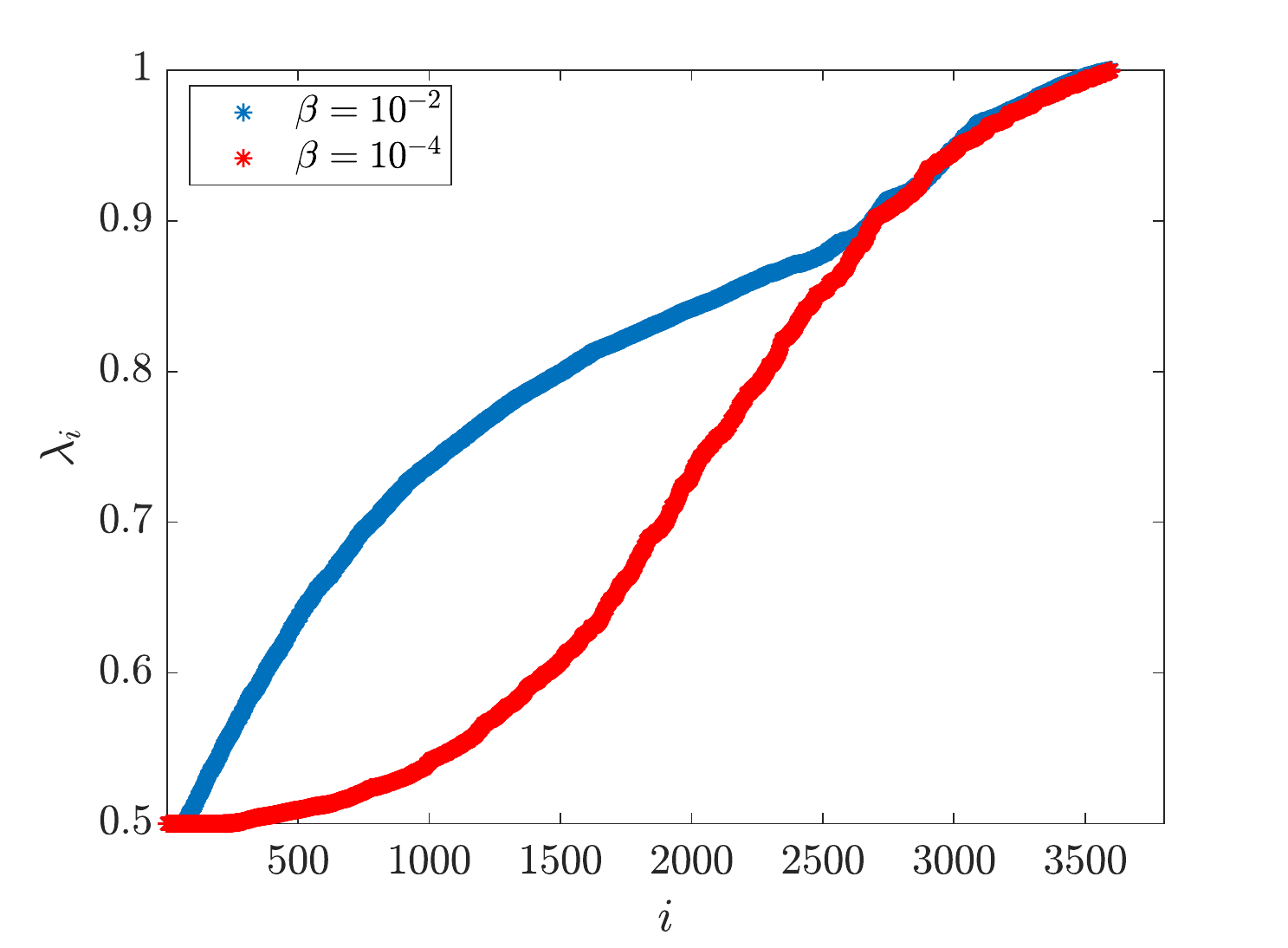}
\includegraphics[width=0.49\linewidth]{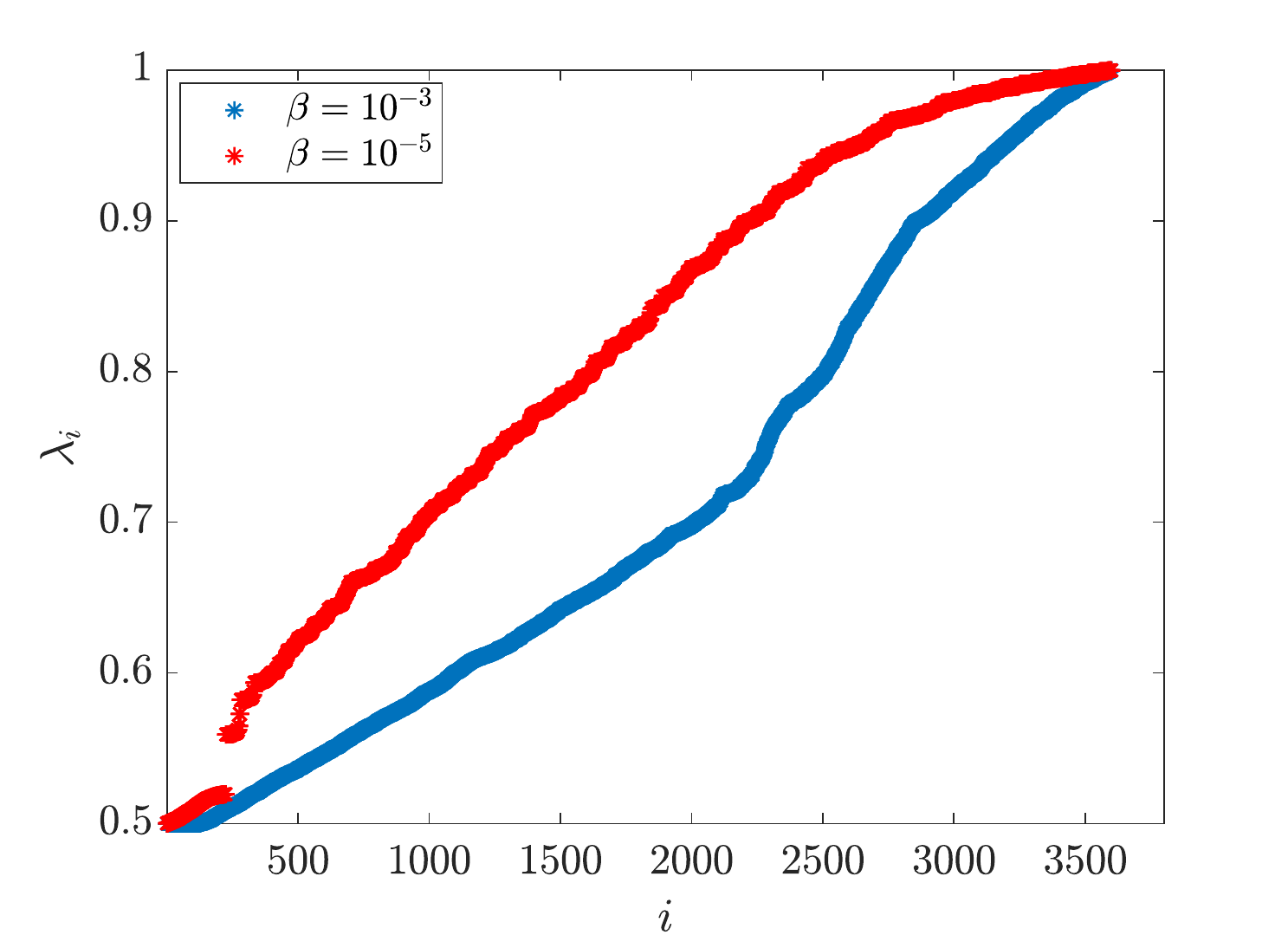}
\end{figure}

In Figure \ref{fig:eigs} we report the eigenvalue distribution of $\widetilde{S}_{\mathrm{int}} ^{-1} S_{\mathrm{int}}$ for a range of values of $\beta$ with diffusion coefficient $\epsilon=\frac{1}{100}$, for a particular Dirichlet test problem.

Further, using Theorem \ref{spectrum_S_int} and \eqref{Rayleigh_quotient}, we can prove the following:
\begin{corollary}\label{spectrum_S}
Let $S$ and $\widehat{S}$ be defined as in \eqref{schur_complement} and \eqref{hat_S_def}, with the matrices defined as in Theorem \ref{spectrum_S_int}, and $T_1$, $T_2$ as in \eqref{linear_transf_T}. Then, given Assumption \ref{assumption_w}:
\begin{displaymath}
\lambda(\widehat{S}^{-1} S) \in \left[ \dfrac{1}{2}, 1\right].
\end{displaymath}
\end{corollary}

By Corollary \ref{spectrum_S}, the matrix $\widehat{S}$ in \eqref{hat_S_def} is an effective approximation of the Schur complement $S$ defined in \eqref{schur_complement}. We highlight that solving (exactly) a system involving the matrix $\widehat{S}$ is costly, so we look for a cheap approximation of the effect of $\widehat{S}^{-1}$ on a generic vector. From \eqref{hat_S_def}, it is clear that the bulk of the work involves approximately applying the inverse of $\Lambda_{21}+ \widehat{M}$ and its transpose. From \eqref{Lambda_21} and \eqref{Mhat}, we have that $\Lambda_{21}+ \widehat{M}$ is block-lower triangular, and a cheap application of its inverse on a vector is given by block-forward substitution, with each block diagonal approximated using a fixed number of V-cycles of a multigrid routine, for example. The matrix $(\Lambda_{21}+ \widehat{M})^\top$ can be handled analogously using block-backward substitution.

\begin{remark}
Theorem \ref{spectrum_S_int} and Corollary \ref{spectrum_S} also hold if no stabilization is applied to the convection--diffusion control problem.
\end{remark}

\begin{remark}
Let us now briefly discuss the applicability of our method if Assumption \ref{assumption_w} does not hold, that is $\mathbf{w}\cdot\mathbf{n}\neq0$ on a portion of $\partial\Omega_N$. As above, this results in the matrix $N$ not being skew-symmetric. Therefore, assuming (non-trivially) a discretization of the adjoint operator has been performed such that one is examining a symmetric matrix system, from \eqref{matrix_H_mass_matrix} we can derive that
\begin{displaymath}
\widetilde{L}= 2 (\epsilon K + W) + H \succeq 2 (\epsilon K + W) - c M_{\partial \Omega_N},
\end{displaymath}
so a sufficient condition for the matrix $\widetilde{L}$ to be positive semi-definite is that $2 (\epsilon K + W) \succeq c M_{\partial \Omega_N}$. Even if the matrix $\widetilde{L}$ were not positive semi-definite, one could also argue that the upper bound on the eigenvalues in Theorem \ref{spectrum_S_int} (and hence that of Theorem \ref{spectrum_S}) would be only slightly larger than $1$ for moderate $\left\|\mathbf{w}\right\|$ and large diffusion coefficient $\epsilon$, as the contribution (in terms of eigenvalues) of the matrix $H$ to $\widetilde{L}$ is no greater than that of the (very sparse) mass matrix on $\partial\Omega_N$ multiplied by $\left\|\mathbf{w}\right\|$.

If the discretization of the adjoint operator is not such that the resulting matrix system is not symmetric, we would need to apply a non-symmetric Krylov solver such as GMRES \cite{Saad_Schultz}, however it is reasonable to suppose that the system will not be highly non-normal. Therefore, we would expect a preconditioner of similar form to that derived in this section could be potent for such problems, but we highlight that the convergence of the Krylov method would not be guaranteed by the clustering of eigenvalues of $\widehat{A}^{-1}A$ and $\widehat{S}^{-1}S$, as for symmetric problems.
\end{remark}

\section{Numerical Results}
\label{5}
We now provide numerical evidence of the effectiveness of our preconditioning strategy. In Section \ref{5_1}, we show how this preconditioner results in more rapid convergence than the state-of-the-art backward Euler solver with the existing preconditioner of \cite{Pearson_Stoll_Wathen}, for the heat control problem. In Section \ref{5_2}, we show the robustness of our strategy with respect to all the parameters involved, for the time-dependent convection--diffusion control problem.

In all our tests we consider only Dirichlet boundary conditions (i.e., $\partial \Omega_N = \emptyset$), but we emphasize again that the method is easily generalized to Neumann and mixed boundary conditions (with the caveats previously outlined for convection--diffusion control). We implement a finite element method, using $\mathbf{Q}_1$ basis functions for state, control, and adjoint variables. As discussed in Section \ref{4_2}, when approximating the $(1,1)$-block it is trivial to invert both matrices $T_1$ and $T_2=T_1^\top$, and we apply 20 steps of Chebyshev semi-iteration to each mass matrix on the diagonal of $A_D$. For the approximation of the Schur complement, we employ block-forward and block-backward substitution to solve for the matrix $\Lambda_{21}+ \widehat{M}$ and its transpose, and approximate each block on the diagonal with 3 V-cycles of the AGMG algebraic multigrid routine \cite{Napov_Notay, Notay1, Notay2, AGMG}. All tests are run on MATLAB R2018b, using a 1.70GHz Intel quad-core i5 processor on an Ubuntu 18.04.1 LTS operating system.

\subsection{Heat control}
\label{5_1}
We first benchmark our method against the backward Euler method coupled with the bespoke, mesh- and $\beta$-independent preconditioner derived in \cite{Pearson_Stoll_Wathen}, for the heat control problem \eqref{heat_control_functional}--\eqref{heat_control_constraints}. Here, $d=2$ (so $\mathbf{x}=[x_1 ,~ x_2 ]^\top$), $\Omega=(-1,1)^2$, $t_{f}=2$, $f=0$, and
\begin{equation*}
\widehat{y}(x_1,x_2,t) = 1+\left[ \left( \dfrac{2}{\pi^2 \beta} + \dfrac{\pi^2}{2} \right) \, e^{t_{f}} + \left( 1 - \dfrac{2}{(2+\pi^2)\beta} - \dfrac{\pi^2}{2} \right) \, e^{t} \right] s(x_1,x_2),
\end{equation*}
where $s(x_1,x_2)= \cos{(\frac{\pi{}x_1}{2})} \cos{(\frac{\pi{}x_2}{2})}$. The analytic solutions to this problem are:
\begin{equation*}
\begin{array}{l}
\vspace{2ex}
y(x_1,x_2,t) = 1+\left( \dfrac{2}{\pi^2 \beta} e^{t_{f}} - \dfrac{2}{(2+\pi^2)\beta}e^{t} \right)s(x_1,x_2),\\
p(x_1,x_2,t) = \left( e^{t_{f}} - e^{t} \right)s(x_1,x_2),
\end{array}
\end{equation*}
with initial condition $y_0$ and Dirichlet boundary condition $g$ obtained from this $y$. We observe that in this case the discretized system obtained after applying Crank--Nicolson is \eqref{modified_system}, with $L=K$.

In our tests, given a (spatial) uniform grid of mesh-size $h=2^{1-l}$, with $l$ the level of grid refinement, we set $\tau = h^2$ for the backward Euler scheme and $\tau = h$ for Crank--Nicolson, motivated by the predicted convergence rates of these methods. We test both schemes with a range of values of $h$ and $\beta$, running preconditioned MINRES to a tolerance of $10^{-6}$ on the relative preconditioned residual norm. We take $\xi=10^{-3}$ in \eqref{prec_Euler}--\eqref{hat_M_Euler} for the backward Euler implementation. Tables \ref{table1}--\ref{table3} present the number of MINRES iterations $\texttt{it}$ required by each method for a range of $\beta$, the CPU time taken in seconds, and the relative errors $y_{\text{error}}$ and $p_{\text{error}}$ (in the scaled vector $\ell^\infty$-norm) obtained for the state and adjoint variables.

\begin{table}[!ht]
\caption{Heat control problem: MINRES iterations, CPU times, and resulting relative errors in $y$ and $p$, for $\beta=10^{-2}$.}\label{table1}
\begin{center}
{\begin{tabular}{|c||c|c|c|c||c|c|c|c|}\hline
\multicolumn{1}{|c||}{} & \multicolumn{4}{c||}{\textbf{Backward Euler}} & \multicolumn{4}{c|}{\textbf{Crank--Nicolson}}\\
\cline{2-9}
$l$  & $\texttt{it}$ & CPU & $y_{\text{error}}$ & $p_{\text{error}}$ & $\texttt{it}$ & CPU & $y_{\text{error}}$ & $p_{\text{error}}$ \\
\hline
\hline
$3$ & 22 & 0.90 & 2.8907e-02 & 9.3751e-03 & 17 & 0.19 & 6.1670e-03 & 9.8782e-03\\
\hline
$4$ & 23 & 4.8 & 1.4472e-03 & 2.3445e-03 & 20 & 0.56 & 1.3137e-03 & 2.2102e-03\\
\hline
$5$ & 26 & 16.2 & 3.6034e-04 & 5.8551e-04 & 20 & 2.5 & 3.2490e-04 & 5.4963e-04 \\
\hline
$6$ & 26 & 340 & 9.0047e-05 & 1.4630e-04 & 22 & 12.2 & 8.1001e-05 & 1.3721e-04 \\
\hline
$7$ & $\dagger$\footnotemark & -- & -- & -- & 24 & 143 & 2.0241e-05 & 3.4423e-05 \\
\hline
\end{tabular}}
\end{center}
\end{table}

\begin{table}[!ht]
\caption{Heat control problem: MINRES iterations, CPU times, and resulting relative errors in $y$ and $p$, for $\beta=10^{-3}$.}\label{table2}
\begin{center}
{\begin{tabular}{|c||c|c|c|c||c|c|c|c|}\hline
\multicolumn{1}{|c||}{} & \multicolumn{4}{c||}{\textbf{Backward Euler}} & \multicolumn{4}{c|}{\textbf{Crank--Nicolson}}\\
\cline{2-9}
$l$  & $\texttt{it}$ & CPU & $y_{\text{error}}$ & $p_{\text{error}}$ & $\texttt{it}$ & CPU & $y_{\text{error}}$ & $p_{\text{error}}$ \\
\hline
\hline
$3$ & 26 & 0.60 & 1.3768e-02 & 1.3732e-02 & 18 & 0.09 & 9.6046e-04 & 1.5371e-02\\
\hline
$4$ & 24 & 8.2 & 2.3058e-03 & 3.4227e-03 & 20 & 0.67 & 1.9959e-04 & 3.5192e-03\\
\hline
$5$ & 25 & 20.4 & 5.0848e-04 & 8.5306e-04 & 21 & 2.5 & 3.7310e-05 & 8.0236e-04 \\
\hline
$6$ & 25 & 327 & 1.2652e-04 & 2.1316e-04 & 25 & 14.3 & 9.2543e-06 & 2.0063e-04\\
\hline
$7$ & $\dagger$ & -- & -- & -- & 28 & 159 & 2.3180e-06 & 5.0100e-05\\
\hline
\end{tabular}}
\end{center}
\end{table}

\begin{table}[!ht]
\caption{Heat control problem: MINRES iterations, CPU times, and resulting relative errors in $y$ and $p$, for $\beta=10^{-4}$.}\label{table3}
\begin{center}
{\begin{tabular}{|c||c|c|c|c||c|c|c|c|}\hline
\multicolumn{1}{|c||}{} & \multicolumn{4}{c||}{\textbf{Backward Euler}} & \multicolumn{4}{c|}{\textbf{Crank--Nicolson}}\\
\cline{2-9}
$l$  & $\texttt{it}$ & CPU & $y_{\text{error}}$ & $p_{\text{error}}$ & $\texttt{it}$ & CPU & $y_{\text{error}}$ & $p_{\text{error}}$ \\
\hline
\hline
$3$ & 17 & 0.40 & 1.0652e-02 & 1.4138e-02 & 16 & 0.08 & 6.2442e-04 & 1.6239e-02 \\
\hline
$4$ & 24 & 8.2 & 1.0722e-03 & 3.5073e-03 & 19 & 0.62 & 2.2555e-04 & 3.5765e-03\\
\hline
$5$ &  26 & 21.0 & 1.7105e-04 & 8.7596e-04 & 22 & 2.74 & 5.7256e-05 & 8.7167e-04\\
\hline
$6$ & 23 & 300 & 4.1337e-05 & 2.1972e-04 & 24 & 13.7 & 1.4088e-05 & 2.1208e-04\\
\hline
$7$ & $\dagger$ & -- & -- & -- & 28 & 158 & 3.7682e-06 & 5.1806e-05 \\
\hline
\end{tabular}}
\end{center}
\end{table}

We see from Tables \ref{table1}--\ref{table3} that the Crank--Nicolson approach with our new preconditioner achieves more accurate solutions than the backward Euler method, in lower CPU time. Comparing the results for grid refinement $l=4,5,6$, that is, for $h=2^{-3},2^{-4},2^{-5}$, this can occur in orders of magnitude lower CPU time for the tests presented here. This is clear because the size of the system required to obtain a fixed accuracy should grow like $\mathcal{O}(h^{-4})$ for backward Euler as opposed to $\mathcal{O}(h^{-3})$ for Crank--Nicolson, and the effectiveness of our new preconditioner allows this to materialize in terms of CPU time. For instance, with this $\Omega$ and $t_{f}$, the choice $h=2^{-5}$ leads to a (Schur complement) system of dimension $8,132,481$ for backward Euler and $254,016$ for Crank--Nicolson; as our preconditioner is optimal for all values of $\beta$, $h$, and $\tau$, this leads to a substantial speed-up. Further, our preconditioned Crank--Nicolson is also able to obtain a solution for levels of refinement (e.g., $l=7$) for which preconditioned backward Euler runs out of memory. We can therefore conclude that the Crank--Nicolson method, coupled with our new preconditioner, is significantly more potent that the widely-used preconditioned backward Euler method for all values of $\beta$.
\footnotetext{$\dagger$ means that the solver ran out of memory.}

\begin{figure}[!htb]
\centering
\caption{Convection--diffusion control problem: Numerical solutions for state and adjoint variables at time $t=1$, with $\epsilon=\frac{1}{20}$, $\beta=10^{-2}$, and $l=5$.}
\label{fig:state_adjoint_1}
\begin{subfigure}[b]{0.45\textwidth}
\centering
\includegraphics[width=0.9\linewidth]{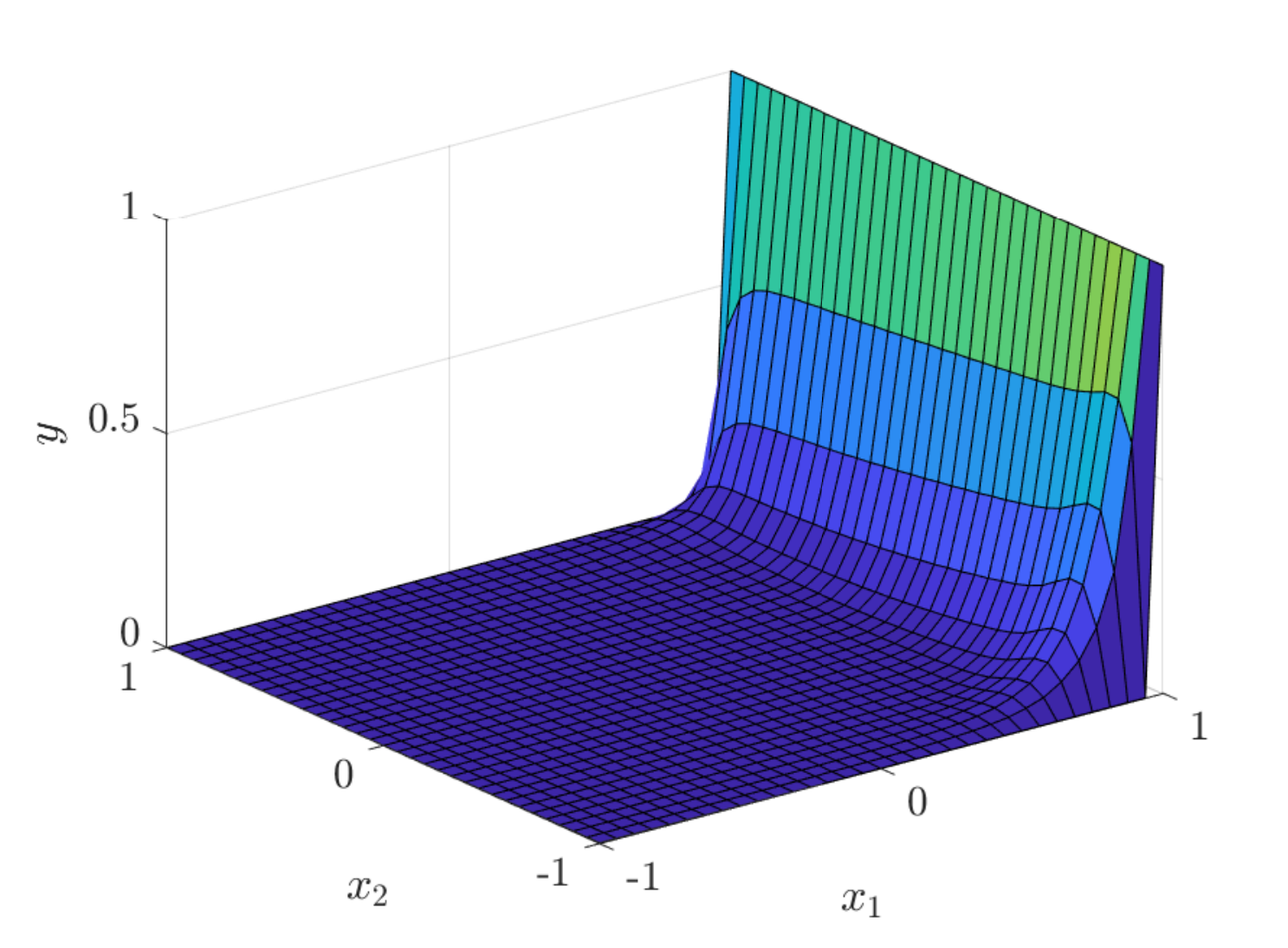}
\caption{State $y$}
\end{subfigure}
\begin{subfigure}[b]{0.45\textwidth}
\centering
\includegraphics[width=0.9\linewidth]{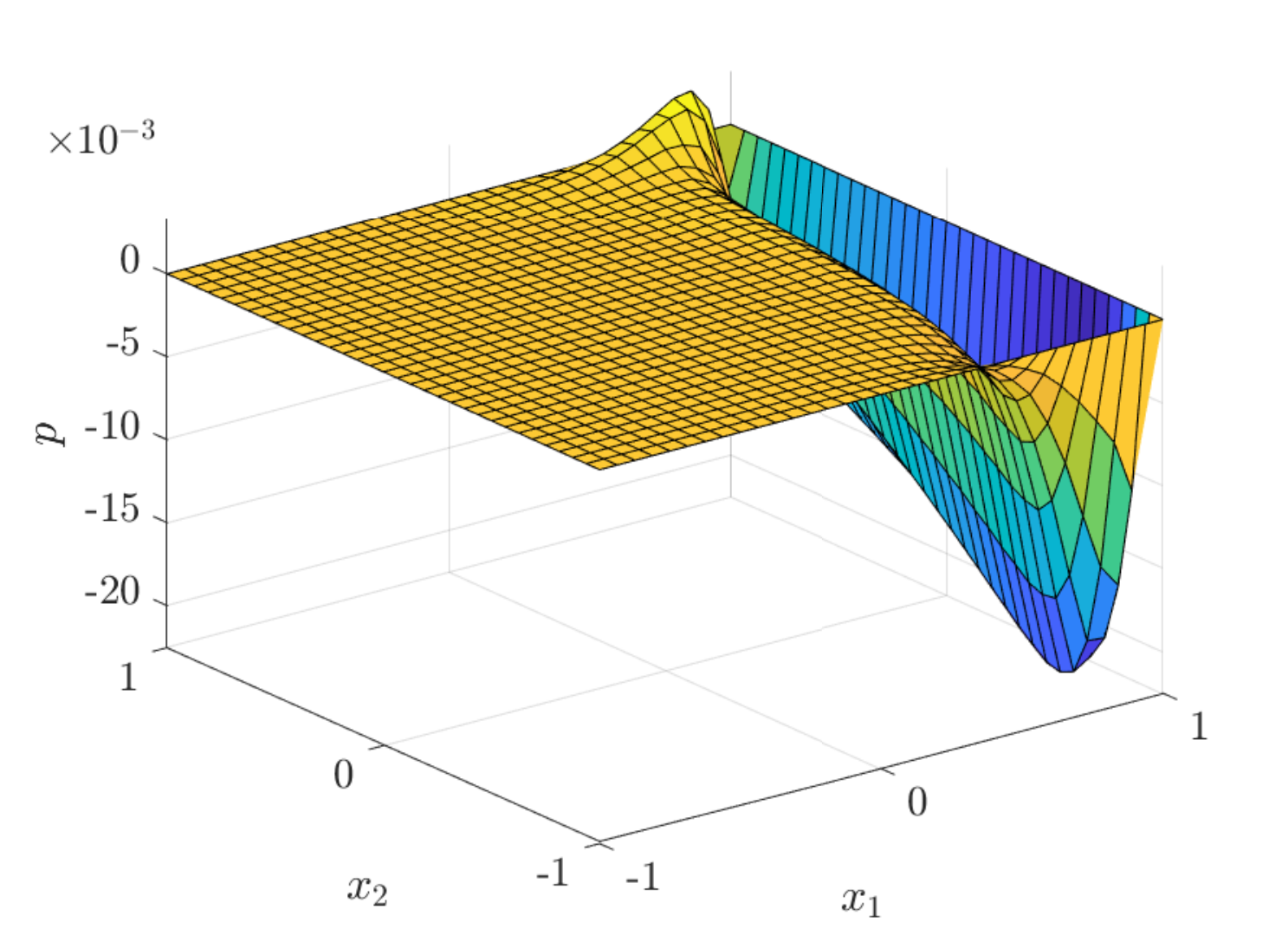}
\caption{Adjoint variable $p$}
\end{subfigure}
\end{figure}

\begin{figure}[!htb]
\centering
\caption{Convection--diffusion control problem: Numerical solutions for state and adjoint variables at time $t=1$, with $\epsilon=\frac{1}{100}$, $\beta=10^{-2}$, and $l=5$.}
\label{fig:state_adjoint_2}
\begin{subfigure}[b]{0.45\textwidth}
\centering
\includegraphics[width=0.9\linewidth]{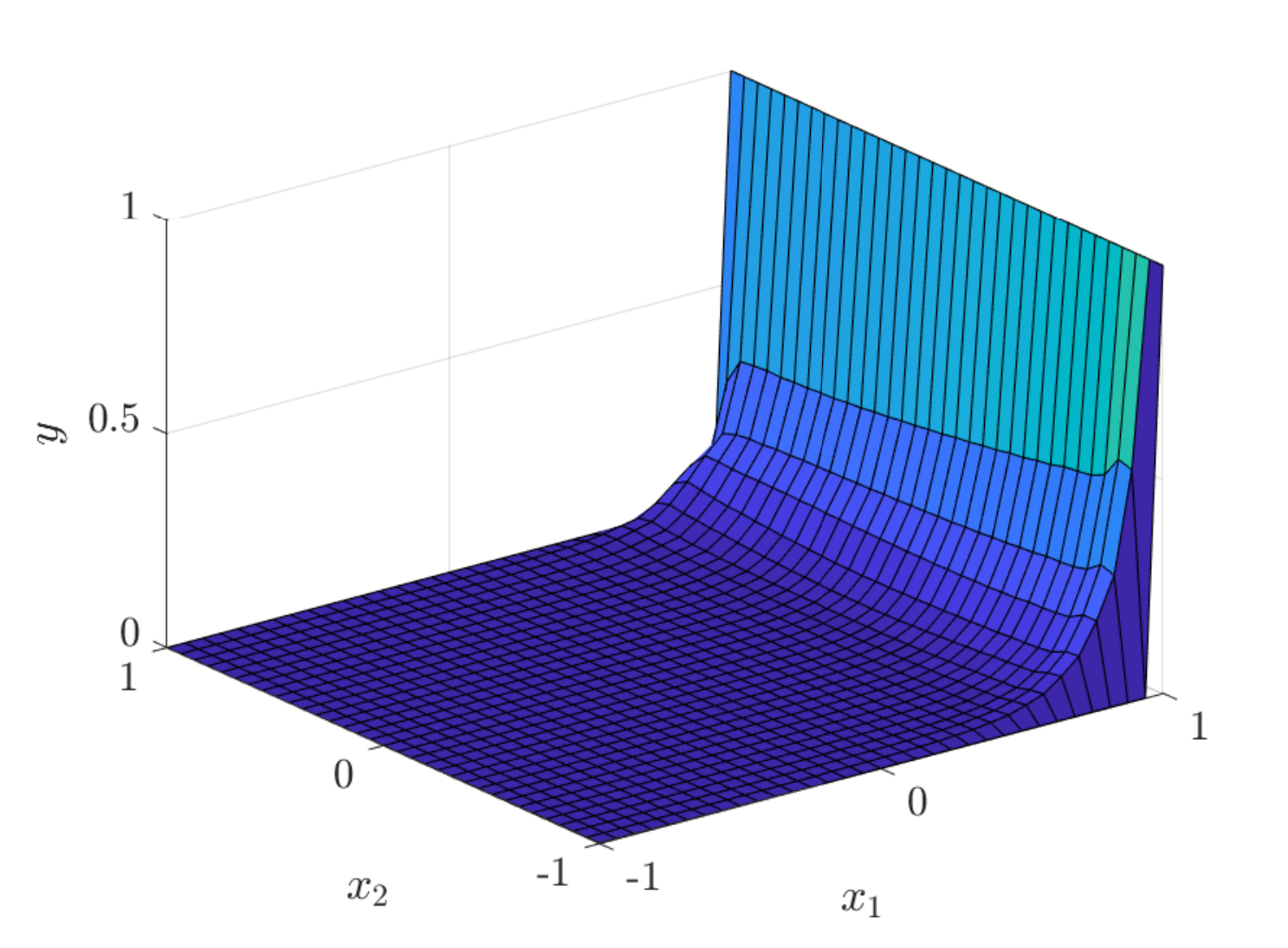}
\caption{State $y$}
\end{subfigure}
\begin{subfigure}[b]{0.45\textwidth}
\centering
\includegraphics[width=0.9\linewidth]{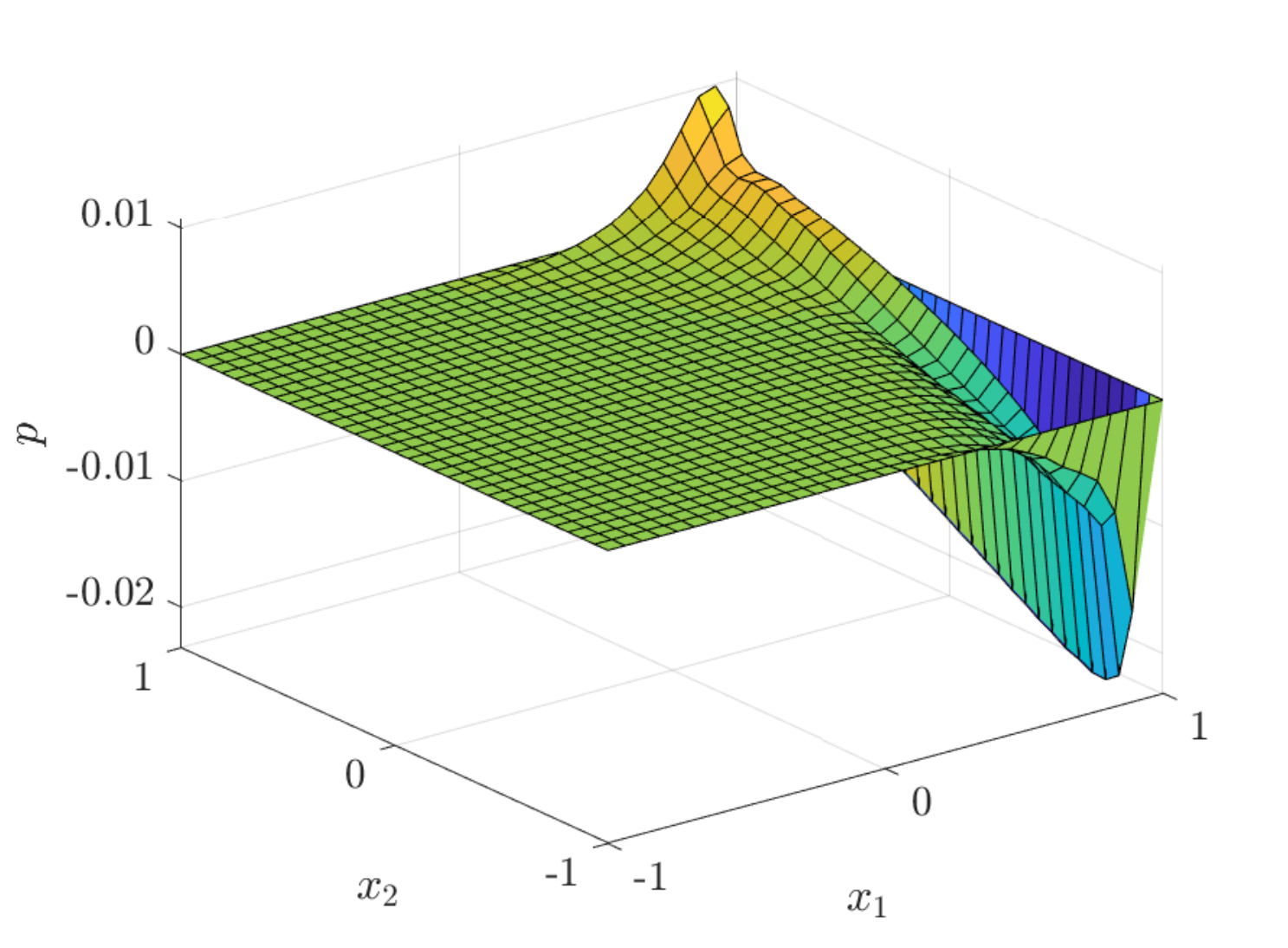}
\caption{Adjoint variable $p$}
\end{subfigure}
\end{figure}

\subsection{Time-dependent convection--diffusion control}
\label{5_2}
We now demonstrate the effectiveness and robustness of our new approach with the respect to the parameters involved in the convection--diffusion control problem \eqref{convection_diffusion_control_functional}--\eqref{convection_diffusion_control_constraints}, specifically $h$ ($=\tau$), $\beta$, $\epsilon$. We again take $d=2$, $\Omega=(-1,1)^2$, $t_{f}=2$, and $f=0$, with wind $\mathbf{w} = [ 2x_2(1-x_1^2),~-2x_1(1-x_2^2) ]^\top$. The initial condition on the state $y$ is given by
\begin{displaymath}
y(x_1,x_2,0)=
\left\{
\begin{array}{ll}
1 & \mathrm{if} \: x_1 = 1,\\
0 & \mathrm{otherwise}.
\end{array}
\right.
\end{displaymath}
Setting $\partial \Omega_1 := \{1\} \times [-1, 1]$ and $\partial \Omega_2=\partial \Omega \setminus \partial \Omega_1$, the boundary condition is given by
\begin{displaymath}
y(x_1,x_2,t)=
\left\{
\begin{array}{ll}
1 & \mathrm{on} \: \partial \Omega_1 \times (0,t_{f}),\\
0 & \mathrm{on} \: \partial \Omega_2 \times (0,t_{f}).
\end{array}
\right.
\end{displaymath}
Finally, the desired state is given by
\begin{displaymath}
\widehat{y} (x_1,x_2,t) = e^{-10(1-x_1)}.
\end{displaymath}
We again run preconditioned MINRES to a tolerance of $10^{-6}$, constructing a (spatial) uniform grid of mesh-size $h=2^{1-l}$ at level $l$, and setting $\tau = h$ in all the tests presented. In Figure \ref{fig:state_adjoint_1}--\ref{fig:state_adjoint_2} we show the numerical solutions for the state $y$ and adjoint variable $p$ at time $t=1$, for $\beta=10^{-2}$ and $l=5$, with both $\epsilon=\frac{1}{20}$ and $\epsilon=\frac{1}{100}$. In Tables \ref{table4}--\ref{table6} we report the number of iterations $\texttt{it}$ required for achieving convergence with values of $\epsilon = \frac{1}{20}, \:\frac{1}{100},\: \mathrm{and} \: \frac{1}{500}$, for a range of $l$ and $\beta$.

\begin{table}[!htb]
\caption{Convection--diffusion control problem: MINRES iterations required with $\epsilon = \frac{1}{20}$, for a range of $l$ and $\beta$.}\label{table4}
\begin{center}
{\begin{tabular}{|c||c|c|c|c|c|c|}\hline
\multicolumn{1}{|c||}{} & \multicolumn{6}{c|}{$~\beta~$}\\
\cline{2-7}
$\displaystyle ~l~$ & $~10^{-1}~$ & $~10^{-2}~$ & $~10^{-3}~$ & $~10^{-4}~$ & $~10^{-5}~$ & $~10^{-6}~$\\
\hline
\hline
$~3~$ & 22 & 24 & 25 & 21 & 15 & 10\\
\hline
$~4~$ & 22 & 25 & 26 & 24 & 19 & 13\\
\hline
$~5~$ & 20 & 25 & 26 & 26 & 23 & 17\\
\hline
$~6~$ & 22 & 23 & 25 & 26 & 26 & 24\\
\hline
$~7~$ & 24 & 24 & 25 & 25 & 26 & 26\\
\hline
\end{tabular}}
\end{center}
\end{table}

\begin{table}[!htb]
\caption{Convection--diffusion control problem: MINRES iterations required with $\epsilon = \frac{1}{100}$, for a range of $l$ and $\beta$.}\label{table5}
\begin{center}
{\begin{tabular}{|c||c|c|c|c|c|c|}\hline
\multicolumn{1}{|c||}{} & \multicolumn{6}{c|}{$~\beta~$}\\
\cline{2-7}
$\displaystyle ~l~$ & $~10^{-1}~$ & $~10^{-2}~$ & $~10^{-3}~$ & $~10^{-4}~$ & $~10^{-5}~$ & $~10^{-6}~$\\
\hline
\hline
$~3~$ & 23 & 26 & 23 & 21 & 15 & 10\\
\hline
$~4~$ & 25 & 26 & 26 & 25 & 21 & 15\\
\hline
$~5~$ & 24 & 27 & 26 & 26 & 23 & 18\\
\hline
$~6~$ & 24 & 25 & 26 & 26 & 26 & 22\\
\hline
$~7~$ & 23 & 25 & 26 & 26 & 26 & 25\\
\hline
\end{tabular}}
\end{center}
\end{table}

\begin{table}[!htb]
\caption{Convection--diffusion control problem: MINRES iterations required with $\epsilon = \frac{1}{500}$, for a range of $l$ and $\beta$.}\label{table6}
\begin{center}
{\begin{tabular}{|c||c|c|c|c|c|c|}\hline
\multicolumn{1}{|c||}{} & \multicolumn{6}{c|}{$~\beta~$}\\
\cline{2-7}
$\displaystyle ~l~$ & $~10^{-1}~$ & $~10^{-2}~$ & $~10^{-3}~$ & $~10^{-4}~$ & $~10^{-5}~$ & $~10^{-6}~$\\
\hline
\hline
$~3~$ & 23 & 26 & 23 & 21 & 15 & 10\\
\hline
$~4~$ & 25 & 26 & 26 & 25 & 21 & 15\\
\hline
$~5~$ & 25 & 27 & 27 & 26 & 25 & 19\\
\hline
$~6~$ & 26 & 27 & 27 & 27 & 27 & 22\\
\hline
$~7~$ & 26 & 27 & 27 & 27 & 25 & 25\\
\hline
\end{tabular}}
\end{center}
\end{table}

As can be seen from Tables \ref{table4}--\ref{table6}, our new preconditioner is highly effective and robust, leading to convergence for all tests in at most $27$ iterations. For $\beta=10^{-5}$ or $10^{-6}$, and larger values $h$, convergence is achieved in a lower number of iterations: this is not surprising as for these values the Schur complement is spectrally `close' to a mass matrix, making the problem easier to solve. Apart from this, we notice that the number of iterations is independent of the parameters involved. We therefore deduce that our method is a potent one for the resolution of time-dependent convection--diffusion control problems, a class of problems which consists of substantial numerical difficulties. The number of iterations required to solve these problems is independent of mesh-size $h$, regularization parameter $\beta$, and diffusion coefficient $\epsilon$.

\section{Conclusions}
\label{6}
We have applied an optimize-then-discretize strategy to tackle the optimal control of time-dependent PDEs, coupled with a Crank--Nicolson scheme in time. We have devised an invertible linear transformation that symmetrizes the resulting matrix system, and derived a new, fast, and robust preconditioner for the saddle-point matrix, which possesses a complex structure. We also have proved that the Schur complement approximation used is optimal with the respect to all parameters involved through bounds on the eigenvalues, and therefore that the preconditioner is optimal and scales linearly in CPU time with respect to matrix dimension. Finally, we have presented numerical results to demonstrate the effectiveness and speed of our preconditioned Crank--Nicolson method. Future work includes extending this new approach to problems with more complex PDE constraints, boundary control problems, and problems with additional algebraic constraints on the state and control variables.

\vspace{1em}
\textbf{Acknowledgements.}~~SL acknowledges financial support from a School of Mathematics PhD studentship at the University of Edinburgh. JWP acknowledges support from the Engineering and Physical Sciences Research Council (EPSRC) grant EP/S027785/1.

%
%
%

\end{document}